\documentclass[onefignum,onetabnum]{siamonline220329}

\usepackage{amsmath,amsbsy,amsfonts,amssymb}
\usepackage{mathtools}

\newcommand{\citeCount}[1]{}
\def\ba#1\ea{\begin{align}#1\end{align}}

\def\bas#1\eas{\begin{align*}#1\end{align*}}

\def\bat#1\eat{\begin{alignat}{3}#1\end{alignat}}

\def\bats#1\eats{\begin{alignat*}{3}#1\end{alignat*}}

\newcommand{\bse}{\begin{subequations}}
\newcommand{\ese}{\end{subequations}}

\newcommand{\nv}{\mathbf{ n}}

\newcommand{\xv}{\mathbf{ x}}
\newcommand{\yv}{\mathbf{ y}}

\newcommand{\grad}{\nabla}

\clearpage

\newlength{\tfwidth}
\newlength{\tfheight}
\newlength{\tfxa}
\newlength{\tfxb}
\newlength{\tfya}
\newlength{\tfyb}
\newcommand{\trimFigWithBox}[6]{%
\setlength\fboxsep{0pt}%
\setlength\fboxrule{1.0pt}
\fbox{\includegraphics[width=#2, clip, trim=#3 #4 #5 #6]{#1}}%
}
\newcommand{\trimFigNoBox}[6]{%
\setlength\fboxsep{1pt}
\setlength\fboxrule{0.0pt}
\fbox{\includegraphics[width=#2, clip, trim=#3 #4 #5 #6]{#1}}%
}
\newcommand{\trimFigHeightWithBox}[6]{%
\setlength\fboxsep{0pt}%
\setlength\fboxrule{1.0pt}
\fbox{\includegraphics[height=#2, clip, trim=#3 #4 #5 #6]{#1}}%
}
\newcommand{\trimFigHeightNoBox}[6]{%
\setlength\fboxsep{1pt}
\setlength\fboxrule{0.0pt}
\fbox{\includegraphics[height=#2, clip, trim=#3 #4 #5 #6]{#1}}%
}

\newsavebox\figBox

\newcommand{\trimw}[6]{%
\sbox\figBox{\includegraphics{#1}}
\setlength{\tfwidth}{\the\wd\figBox}
\setlength{\tfheight}{\the\ht\figBox}
\setlength{\tfxa}{\tfwidth*\real{#3}}%
\setlength{\tfxb}{\tfwidth*\real{#4}}%
\setlength{\tfya}{\tfheight*\real{#5}}%
\setlength{\tfyb}{\tfheight*\real{#6}}%
\trimFigNoBox{#1}{#2}{\tfxa}{\tfya}{\tfxb}{\tfyb}%
}

\newcommand{\trimwb}[6]{%

\sbox\figBox{\includegraphics{#1}}
\setlength{\tfwidth}{\the\wd\figBox}
\setlength{\tfheight}{\the\ht\figBox}
\setlength{\tfxa}{\tfwidth*\real{#3}}%
\setlength{\tfxb}{\tfwidth*\real{#4}}%
\setlength{\tfya}{\tfheight*\real{#5}}%
\setlength{\tfyb}{\tfheight*\real{#6}}%
\trimFigWithBox{#1}{#2}{\tfxa}{\tfya}{\tfxb}{\tfyb}%
}

\newcommand{\trimh}[6]{%
\sbox\figBox{\includegraphics{#1}}
\setlength{\tfwidth}{\the\wd\figBox}
\setlength{\tfheight}{\the\ht\figBox}
\setlength{\tfxa}{\tfwidth*\real{#3}}%
\setlength{\tfxb}{\tfwidth*\real{#4}}%
\setlength{\tfya}{\tfheight*\real{#5}}%
\setlength{\tfyb}{\tfheight*\real{#6}}%
\trimFigHeightNoBox{#1}{#2}{\tfxa}{\tfya}{\tfxb}{\tfyb}%
}

\newcommand{\trimhb}[6]{%

\sbox\figBox{\includegraphics{#1}}
\setlength{\tfwidth}{\the\wd\figBox}
\setlength{\tfheight}{\the\ht\figBox}
\setlength{\tfxa}{\tfwidth*\real{#3}}%
\setlength{\tfxb}{\tfwidth*\real{#4}}%
\setlength{\tfya}{\tfheight*\real{#5}}%
\setlength{\tfyb}{\tfheight*\real{#6}}%
\trimFigHeightWithBox{#1}{#2}{\tfxa}{\tfya}{\tfxb}{\tfyb}%
}
\newtheorem{problem}{Problem}
\newtheorem{approach}{Approach}
\usepackage{tikz}

\newcommand{\das}[1]{#1}
\newcommand{\dasb}[1]{#1}

\usetikzlibrary{shapes,arrows}

\newsiamremark{remark}{Remark}

\headers{Newton-based Grad Shafranov}{D. A. Serino, Q. Tang, et al.}

\title{An adaptive Newton-based free-boundary Grad--Shafranov solver\thanks{Submitted to the editors July 2024.
\funding{
This work was jointly supported by the U.S. Department of Energy through the Fusion Theory Program of
the Office of Fusion Energy Sciences and the SciDAC partnership on Tokamak Disruption Simulation between the
Office of Fusion Energy Sciences and the Office of Advanced Scientific Computing. It was also partially supported
by Mathematical Multifaceted Integrated Capability Center (MMICC) of Advanced Scientific Computing Research.
Los Alamos National Laboratory is operated by Triad National Security, LLC, for the National Nuclear Security
Administration of U.S. Department of Energy (Contract No.~89233218CNA000001).
}}}
\author{Daniel A. Serino\thanks{Los Alamos National Laboratory, Los Alamos, NM
  (\email{dserino@lanl.gov},
  \email{qtang@lanl.gov},
  \email{xtang@lanl.gov},
  \email{lipnikov@lanl.gov})}
\and Qi Tang\footnotemark[2] \thanks{\das{Current address: School of Computational Science and Engineering, Georgia Institute of Technology, Atlanta, GA (\email{qtang@gatech.edu})}}
\and Xian-Zhu Tang\footnotemark[2]
\and Tzanio V.~Kolev\thanks{Lawrence Livermore National Laboratory, Livermore, CA
  (\email{kolev1@llnl.gov})}
\and Konstantin Lipnikov\footnotemark[2]
}

\ifpdf
\hypersetup{
  pdftitle={Newton Grad Shafranov},
  pdfauthor={D. A. Serino, Q. Tang, et al.}
}
\fi

\begin{document}

\maketitle

\begin{abstract}
\das{Equilibria} in magnetic confinement devices result from force
balancing between the Lorentz force and the plasma pressure
gradient. In an axisymmetric configuration like a tokamak, such an
equilibrium is described by an elliptic equation for the poloidal
magnetic flux, commonly known as the Grad--Shafranov equation. It is
challenging to develop a scalable and accurate free-boundary
Grad--Shafranov solver, since it is a fully nonlinear optimization
problem that \das{simultaneously} solves for the magnetic field coil current
outside the plasma to control the plasma shape. In this work, we
develop a Newton-based free-boundary Grad--Shafranov solver using
adaptive finite elements and preconditioning strategies. The
free-boundary interaction leads to the evaluation of a
domain-dependent nonlinear form of which its contribution to the
Jacobian matrix is achieved through shape calculus. The optimization
problem aims to minimize the distance between the plasma boundary and
specified control points while satisfying two non-trivial constraints,
which correspond to the nonlinear finite element discretization of the
Grad--Shafranov equation and a constraint on the total plasma current
involving a nonlocal coupling term. The linear system is solved by a
block factorization, and AMG is called for subblock elliptic
operators. The unique contributions of this work include the treatment
of a global constraint, preconditioning strategies, nonlocal
reformulation, and the implementation of adaptive finite elements. It
is found that the resulting Newton solver is robust,
successfully
reducing the nonlinear residual to 1e-6 and lower in a small handful of
iterations
while 
addressing the challenging case to find a Taylor state equilibrium
where conventional Picard-based solvers fail to converge.
\end{abstract}



\section{Introduction}

A high-quality \das{magnetohydrodynamics} (MHD) equilibrium plays a critical role in the whole
device modeling of tokamaks, for both machine design/optimization and
the discharge scenario modeling.  It is also the starting point for
linear stability analyses as well as the time-dependent MHD
simulations to understand the evolution of these instabilities
\cite{bonilla2023fully, jorti2023mimetic}.  The most useful MHD
equilibrium in practice is solved using the well-known free-boundary
Grad--Shafranov equation \cite{jardin2010computational}, which also
determines the coil current in poloidal magnetic field coils to
achieve desired plasma shape.  The conventional approach for this
problem is to use a Picard-based iterative method
\cite{jardin2010computational, jeon2015development}.  For instance, we
developed a cut-cell Grad--Shafranov solver based on Picard iteration
and Aitken's acceleration in \cite{liu2021parallel}; another popular example
is the open-source finite-difference-based solver in
Ref.~\cite{freegsgithub}.
The solver in \cite{liu2021parallel} was
able to provide a less challenging, low-$\beta$ equilibrium for the
whole device simulations in \cite{bonilla2023fully}.  However, we
found that it is difficult to find a Taylor state equilibrium (will be
defined later), as needed in \cite{jorti2023mimetic}, where the plasma
is assumed to be 0-$\beta$.  This work addresses such a challenging
case through developing a Newton-based approach for an optimization
problem to seek a 0-$\beta$ equilibrium.

For challenging cases, it is well known that a Picard-based solver may
have trouble to converge.  In the previous study of
\cite{liu2021parallel}, we found that the relative residual error for
the Taylor state equilibrium often stops reducing in a level of 1e-1
or 1e-2.  The challenge to extend to a Newton-based method is to
compute the Jacobian with the shape of the domain being implicitly
defined by the current solution.  The breakthrough of the Newton-based
solver was recently developed in \cite{heumann2015quasi}, in which a
shape calculus \cite{delfour2011shapes} was involved to compute its
analytical Jacobian.  However, a direct solver was used to address the
linearized system in their code. Although it works fine for a
small-scale 2D problem, it is known that a direct solver is not
scalable in parallel if we seek a more complicated shape control case
like that in the SPARC tokamak \cite{creely-sparc-jpp-2020}.
Ref.~\cite{heumann2015quasi} was later extended to a more advanced
control problem in \cite{blum2019automating}.  More recently,
\das{Ref.~\cite{freegsgithub}} has been extended to a Jacobian-free
Newton--Krylov (JFNK) algorithm as an iterative solver in
\cite{freegsnke}, \das{particularly addressing a dynamic plasma equilibrium problem}.  
However, this did not resolve the
issue, as a JFNK algorithm without a proper preconditioner is known to
be challenging to converge.  Aiming for a case like the Taylor state
equilibrium, the primary goal of this work is to identify a good
preconditioning strategy for the linearized system.  We stress that
the studied cases herein introduce a more difficult optimization
problem than the cases addressed in \cite{heumann2015quasi},
necessitating the preconditioning strategy.

This work is partially inspired by our previous work
\cite{tang2022adaptive} in which we use MFEM \cite{mfem-web,
  anderson2021mfem} to develop an adaptive, scalable, fully implicit
MHD solver.  The flexibility of the solver interface of MFEM was
found to be ideal for developing a complicated nested preconditioning
strategy.  The added-mesh-refinement (AMR) interface provided by MFEM is another attractive
feature for practical problems.  Besides resolving the local features
in a solution, we have seen that a locally refined grid will make the
overall solver more scalable than that on a uniformly refined grid
\cite{tang2022adaptive}.  Both points
have been \das{extensively} explored
during the implementation of this work, aiming for a practical,
scalable solver for more challenging \das{equilibria} with complicated
constraints.  Note that we have developed an adaptive Grad--Shafranov
solver in \cite{dpgGS} but only for a fixed-boundary problem, which is
significantly easier than a free-boundary problem herein.

The contribution of work in this paper can be summarized as follows:
\begin{itemize}
\item We develop a Newton-based free-boundary Grad--Shafranov solver
  with a preconditioned iterative method for the linearized system.
  To the best of our knowledge, the preconditioning strategy has not
  been explored for the Newton-based free-boundary Grad--Shafranov solver.  The
  preconditioning strategy is discovered through several attempts of
  choosing different factorizations and different choices of sub-block
  preconditioners.  The resulting algorithm
  \das{can be further parallelized without major modifications.}
  \das{This is critical for extending the work to more complicated problems}.
\item We propose an alternative objective function for plasma shape
  control which encourages the contour line defining the plasma
  boundary to align closely with specified control points. This
  results in a more challenging nonlinear system to solve in which we
  solve effectively using a quasi-Newton approach.
\item We extend the solver to seeking the Taylor state equilibrium
  that requires to address two non-trivial constraint equations, which
  enforce the total current in the plasma domain and constrain the
  shape of the plasma with a nonlinear condition.  Both constraints
  help to find the optimal equilibrium as the initial condition for
  practical tokamak MHD simulations.  This and the above point
  emphasize the need of a good preconditioner.
\item We develop the proposed algorithm in MFEM, taking the full
  advantage of its scalable finite element implementation, scalable
  sub-block preconditioners, AMR interface, along with many other
  features.
\end{itemize}

\smallskip

The remainder of the manuscript is organized as follows.  In
Section~\ref{sec:gov}, we introduce the Grad--Shafranov equations and
present three prototypical models for plasma current and \das{pressure} profile.
Section~\ref{sec:disc} introduces the finite element discretization
and its linearization for use in a Newton method.  In
Section~\ref{sec:opt}, we formulate the plasma shape control problem
as a constrained nonlinear optimization problem and introduce a Newton
iteration to find solutions.  This reduces the problem to solving a
reduced linear system, for which we discuss effective preconditioning
techniques in Section~\ref{sec:prec}.  Some details of our
implementation based on MFEM are given in Section~\ref{sec:mfem}.  We
demonstrate the performance of the preconditioners, the effectiveness
of the new nonlinear objective function, and the behavior of the AMR
algorithm in Section~\ref{sec:results}.  In
Section~\ref{sec:conclusions}, we conclude and discuss the limitations
of our approach.

\bigskip

\section{Governing Equations}
\label{sec:gov}


The equilibrium of a plasma can be described using the magnetohydrodynamics (MHD) equations,
\bse
\label{eq:mhd}
\begin{align}
  J \times B &= \grad p, \\
  \mu J &= \grad \times B, \label{eq:mhdJ}
\end{align}
\ese
where $J$ is the electric current, $B$ is the magnetic field,
$p$ is the plasma pressure, and $\mu$ is the (constant) magnetic permeability.
We consider an axisymmetric solution to~\eqref{eq:mhd} in \das{cylindrical} coordinates $(r, \phi, z)$.
In an axisymmetric device, the magnetic field is assumed to satisfy
\begin{align}
  B = \frac{1}{r}\grad \psi \times e_\phi + \frac{f(\psi)}{r} e_\phi,
  \label{eq:Bfield}
\end{align}
where $\psi=\psi(r, z)$ is the poloidal magnetic flux, $f$ is the diamagnetic function (which defines
the toroidal magnetic field, $B_\phi = \frac{f(\psi)}{r}$),
and $e_\phi$ is the unit vector for the coordinate $\phi$.
Due to this, the pressure is constant along field lines (i.e., $p=p(\psi)$).
The Grad--Shafranov equation then results from substituting~\eqref{eq:Bfield} into~\eqref{eq:mhd}.
\begin{align}
  \Delta^* \psi:= r \partial_r \left(\frac{1}{r} \partial_r \psi\right) + \partial_z^2 \psi =
  -\mu r^2 p'(\psi) - f(\psi) f'(\psi)
\end{align}

\begin{figure}[h]
  \centering
    \includegraphics[width=\textwidth]{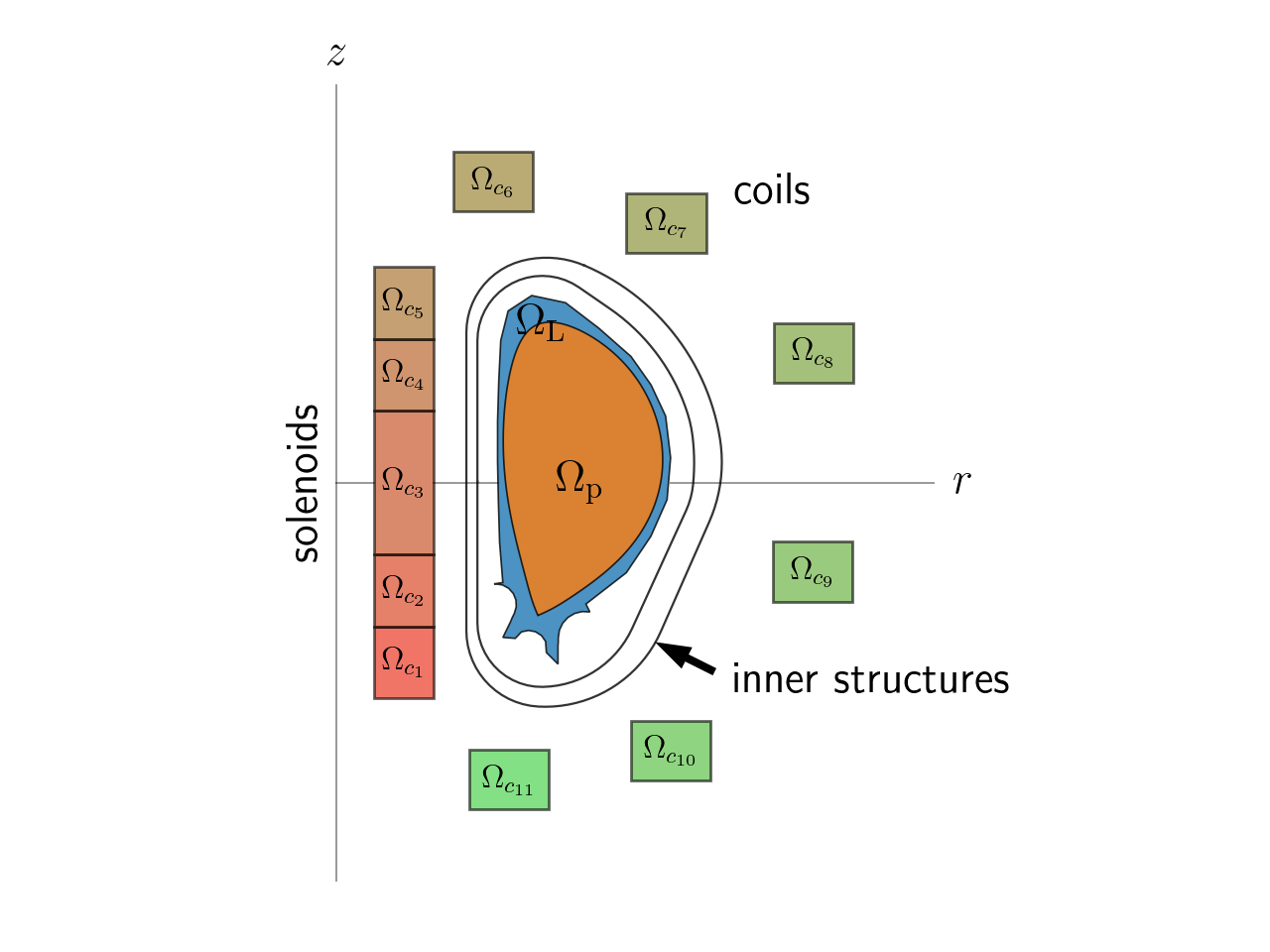}
    \caption{Left: geometry of the problem depicting the $(r, z)$ cross section of a tokamak.
      The domain contains coils and solenoids, $\Omega_{c_i}$, for $i=1,\dots, N_I$,
      a limiter region, $\Omega_{\rm L}$, inside of the first wall,
      and
      the plasma region, $\Omega_{p}$, enclosed by $\Omega_{\rm L}$.
    }
    \label{fig:plasmaDomain}
\end{figure}

Consider a poloidal $(r, z)$ cross section of a tokamak, depicted in Figure~\ref{fig:plasmaDomain}.
Define $\Omega = [0, \infty] \times [-\infty, \infty]$ to be the domain centered with the tokamak geometry.
The tokamak is composed of coils, solenoids and inner structures.
The coils and solenoids are represented by the subdomains
$\Omega_{{\rm C}_i}$ labeled by $i=1,\dots, N_I$.
Let $\Omega_{\rm L}$ be the domain that is enclosed by the first wall.
Each of the coils and solenoids has a current, $I_i$ for $i=1,\dots, N_I$.
\das{We consider an optimization problem where the currents are tuned to produce an equilibrium
  where the boundary of the plasma domain matches a desired reference curve, $\Gamma_p$.
  If the plasma domain is defined by $\Omega_{\rm p}(\psi)\subset\Omega_{\rm L}$,
  then we define an objective function
  $d(\partial \Omega_{\rm p}(\psi), \Gamma_p)$, which represents the ``distance'' between the plasma boundary,
  $\partial \Omega_{\rm p}(\psi)$, and the desired plasma boundary.
  This distance is measured at some predetermined control points. 
  We provide further details on our choice for the objective function in Section~\ref{sec:opt}.}
Additionally, there is a scaling parameter, $\alpha$, involved in the definitions of $f$ and $p$
that is \das{included in the optimization} to \das{constrain} the total plasma current, $I_{\rm p}$.
The total plasma current is represented by a domain integral over $\Omega_{\rm p}(\psi)$.
This control problem, when combined with boundary conditions for the flux at the axis and far-field,
can be summarized as the following PDE-constrained optimization problem.

{\problem Determine the poloidal magnetic flux, $\psi$, and solenoid and coil currents, $I_i$,
  and scaling parameter, $\alpha$,
  that minimizes $d(\partial \Omega_{\rm p}(\psi), \Gamma_p)$
  subject to the constraints
  \bse
  \label{eq:govEq}
  \begin{align}
    - \frac{1}{\mu r} \Delta^* \psi
    &= \begin{cases}
         r p'(\psi) + \frac{1}{\mu r} f(\psi) f'(\psi), & \text{\rm in } \Omega_{\rm p}(\psi), \\
         I_i / | \Omega_{c_i} |, & \text{\rm in } \Omega_{c_i}, \\
         0, & \text{\rm elsewhere in } \Omega_{\infty}
       \end{cases} \\
    \psi(0, z) &= 0, \\
    \lim_{\| (r, z)\| \rightarrow +\infty} \psi(r, z) &= 0 \\
    I_{\rm p} &= \int_{\Omega_{\rm p}(\psi)} \left(r p'(\psi) + \frac{1}{\mu r} {ff'}(\psi)\right) \; dr dz. \label{eq:Ip}
  \end{align}
  \ese
}

The plasma domain, $\Omega_{\rm p}$,  is implicitly determined by $\psi$.
Let $\psi_{\rm ma}$ and $\psi_{\rm x}$ correspond to the $\psi$ values at the
magnetic axis $(r_{\rm ma}(\psi), z_{\rm ma}(\psi))$
and x-point $( r_{\rm x}(\psi), z_{\rm x}(\psi))$, respectively.
The magnetic axis corresponds to the global minimum of $\psi$ inside $\Omega_L$.
When it exists, the x-point is taken to be one of the saddle points of $\psi$ inside the limiter region that has
a $\psi$ value closest to $\psi_{\rm ma}$.
When there are no candidate saddle points inside $\Omega_L$, then this point is chosen to correspond to the maximum
value of $\psi$ on $\partial \Omega_L$.
The level-set of $\psi_{\rm x}$ defines the boundary of the plasma domain.

In order to close the problem, the functions $p$ and $f$ must be supplied. We consider three approaches for defining these functions.
The first option is a simple model used in~\cite{heumann2015quasi} based on Ref.~\cite{luxonbrown}:

{\approach
  \textbf{Luxon and Brown Model}: the functions $p'$ and $ff'$ are provided as
  \begin{align}
    p'(\psi) = \alpha \frac{\beta}{r_0} (1 - \psi_N(\psi)^{\delta})^{\gamma},
    \qquad
    f(\psi) f'(\psi) = \alpha (1 - \beta) \mu r_0 (1 - \psi_N(\psi)^{\delta})^{\gamma},
    \label{eq:lb}
  \end{align}
  where $\psi_{N}(\psi)$ represents a normalized poloidal flux, given by
  \begin{align}
    \psi_N(\psi) = \frac{\psi - \psi_{\rm ma}(\psi)}{\psi_{\rm x}(\psi) - \psi_{\rm ma }(\psi)}.
    \label{eq:psiN}
  \end{align}
Here  $\alpha$ is a scaling constant that is used to satisfy the constraint~\eqref{eq:Ip}. We let $r_0=6.2,\, \delta=2,\, \beta=0.5978,\, \gamma=1.395$,
  using the same coefficients used in \cite{heumann2015quasi} for the ITER configuration.
  This problem is referred to as ``inverse static, with given plasma current $I_{\rm p}$'' in \cite{heumann2015quasi}.
}

\smallskip

 Next we consider the Taylor state equilibrium~\das{\cite{taylorstate}}.
In a tokamak disruption, if the thermal quench is driven by parallel plasma transport
along open magnetic field lines, the plasma current profile is supposed to relax as the result
of magnetic reconnection and magnetic helicity conservation.
Taylor state corresponds to the extreme case in which
$\grad \times B = \alpha B$ with $\alpha$ a global constant.
Substituting the magnetic field representation \eqref{eq:Bfield} and the current in MHD \eqref{eq:mhdJ} into
the Taylor state constraint gives
\begin{align}
 \frac{\partial f}{\partial\psi} \nabla\psi \times \frac{1}{r} e_\phi
- \frac{1}{r} \left(\Delta^*\psi\right)  e_\phi
=
\alpha\left[\frac{1}{r}\grad \psi \times e_\phi + \frac{f(\psi)}{r} e_\phi\right].
\end{align}
Projecting onto two components, one has
\begin{align}
  \frac{\partial f}{\partial\psi} & = \alpha, \label{eq:taylor-dIdpsi} \\
  -\Delta^*\psi & = \alpha  f(\psi). \label{eq:taylor-rjphi}
\end{align}
In the free-boundary Grad--Shafranov equilibrium formulation, 
$f$ is usually a function of the normalized flux $\psi_N(\psi)$.
Solving \eqref{eq:taylor-dIdpsi}, we have
\begin{align}
f(\psi) & = f_{\rm x} - \alpha \left(\psi_{\rm x} - \psi_{\rm ma}\right) + \alpha \left(\psi_{\rm x} - \psi_{\rm ma}\right) {\psi}_N(\psi) \nonumber\\
& =  f_{\rm x} + \alpha (\psi-\psi_{\rm x})
\end{align}
It is interesting to note that
\begin{align}
f(\psi_{\rm x}) = f \big{|}_{{\psi_N(\psi)}=1} = f_{\rm x}
\end{align}
is set by the vacuum toroidal field (given by the current in the toroidal field coils).
The toroidal field on the magnetic axis has an explicit dependence on $\alpha,$ via
\begin{align}
f(\psi_{\rm ma}) = f \big{|}_{{\psi_N(\psi)}=0} = f_{\rm x} + \alpha \left(\psi_{\rm ma} -  \psi_{\rm x} \right).
\end{align}


On ITER, there is no change in $f_{\rm x}$ over the time period of thermal quench and perhaps current quench as well,
so $f_{\rm x}$ is set by the original equilibrium. The factor $\alpha$ sets the total plasma current,
\begin{align}
I_{\rm p}(\alpha) = \int_{\Omega_{\rm p}(\psi)}  \frac{1}{\mu r} {ff'}(\psi)  \; dr dz
\end{align}
\das{In other words, we include $\alpha$ as a free variable in the optimization
  in order to satisfy the constraint
  of matching the given total toroidal plasma current $I_{\rm p}$.
}

In summary, the second option we consider is:
{\approach
  \textbf{Taylor State Equilibrium}:
  $p'(\psi)=0$ and the diamagnetic function is given by
\begin{align}
  f(\psi) = f_{\rm x} + \alpha (\psi - \psi_{\rm x}),
\end{align}
where $f_{\rm x}$ is a constant set by the vacuum toroidal field.
}

\smallskip

\noindent The final option we consider involves  some \das{synthetic} measurement data from a proposed 15MA ITER baseline \cite{liu2015modelling}, which relates $p$ and $f$ to $\psi$.


{\approach
\textbf{15MA ITER baseline case}: $p'$ and $ff'$ are provided using tabular data as a
function of a normalized poloidal flux.
Splines, $S_{p'}(\psi_N)$ and $S_f(\psi_N)$, are fitted to the tabular data and we define the functions as
\begin{align}
  p(\psi) = S_{p'}(\psi_N(\psi))
  \qquad
  f(\psi) = f_{\rm x} + \alpha S_{f}(\psi_N(\psi)),
\end{align}
where $f_{\rm x}$ is the toroidal field function the separatrix, which is constrained by the total poloidal
current in the toroidal field coils outside the chamber wall,
and $S_{f}(\psi_N(\psi))$ is normalized and centered so that $S_{f}(1) = 0$ and $S_{f}(0) = 1$.
Again, $\alpha$ is a scaling constant that is used to satisfy the constraint~\eqref{eq:Ip}
}

\section{Discretization}
\label{sec:disc}


\begin{figure}[tbh]
  \centering
  \includegraphics[width=.8\textwidth]{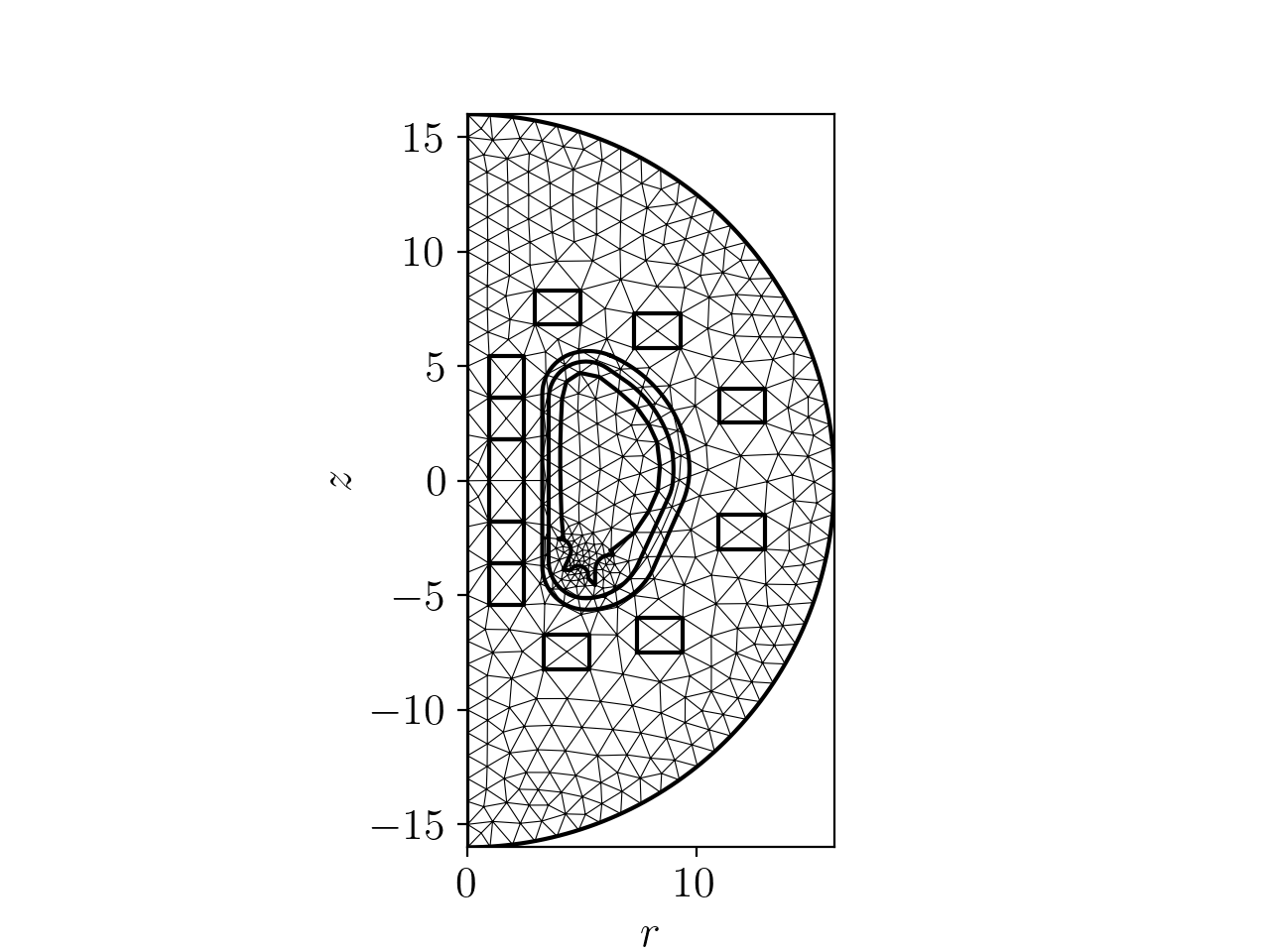}
  \caption{Computational mesh over the region $\Omega$. The geometry of the coils, soilenoids, and inner structures are overlaid over the mesh.\label{fig:initialmesh}}
\end{figure}

Let $\Omega$ denote the computational domain, represented by the semi-circle
$\{r, z: r^2+z^2\le R^2, \; r\ge 0\}$.
Figure~\ref{fig:initialmesh} shows the computational domain and its exemplar coarse mesh.
The domain is meshed using 2D elements (i.e., triangles)
and the solution, $\psi(r,z)$, is assumed to be in $H^1(\Omega)$.
Let $\Gamma_R$ correspond to the arc of the semi-circle, $\{r, z: r^2 + z^2 = R^2, \; r\ge 0 \}$.
Additionally, let $v_i(r,z)$, $i=1,\dots,N_{\rm dof}$ represent the $i^{\rm th}$ shape function.
We wish to solve the following nonlinear equations for $\psi$.
\begin{align}
  {\rm a}(\psi, v_i) = l(I, v_i), \qquad i=1,\dots, N_{\rm dof} . \label{eq:nonlin}
\end{align}
Here, $a$ represents the weak form of~\eqref{eq:govEq}, and
contains contributions from the $\Delta^*$ operator, plasma terms, and far-field boundary terms.
\begin{align}
  a(\psi, v) = a_{\Delta^*}(\psi, v) + a_{\rm plasma}(\psi, v)  + a_{\rm far-field}(\psi, v)
\end{align}
The terms are given by
\bse
\begin{align}
  a_{\Delta^*}(\psi, v) =\;& {\int_{\Omega} \frac{1}{\mu r} \grad \psi \cdot \grad v \; dr dz}\label{eq:terma}\\
  a_{\rm plasma}(\psi, v) =\;&- {\int_{\Omega_{\rm p}(\psi)} \left( r p'(\psi) + \frac{1}{\mu r} f(\psi) f'(\psi) \right) v \; dr dz}\label{eq:termb}\\
  a_{\rm far{\rm -}field}(\psi, v) =\;& \frac{1}{\mu} \int_{\Gamma_R} \psi(\xv) N(\xv) v(\xv) \; dS(\xv) \label{eq:termc}\\
                         &+\frac{1}{2\mu}\int_{\Gamma_R} \int_{\Gamma_R} (\psi(\xv) - \psi(\yv))M(\xv, \yv) (v(\xv) - v(\yv)) \; dS(\xv) dS(\yv),\nonumber
\end{align}
\ese
The far-field terms result from applying appropriate Green's functions on $\Gamma_R$~\cite{blum86} and result in a dense boundary block. These terms use the following definitions
\begin{align}
  M(\xv, \yv) &= \frac{k_{\xv, \yv}}{2 \pi (\xv_r \yv_r)^{3/2}}
                \left(\frac{2 - k^2_{\xv, \yv}}{2 - 2 k^2_{\xv, \yv}} E(k_{\xv, \yv}) - K(k_{\xv, \yv}) \right), \\
  N(\xv) &= \frac{1}{\xv_r} \left(\frac{1}{\delta_+} + \frac{1}{\delta_-} - \frac{1}{R} \right),
\end{align}
where $K$ and $E$ are the complete elliptic integrals of the first and second kind, respectively, and
\begin{align*}
  \delta_{\pm} = \sqrt{\xv_r^2 + (R \pm \xv_z)^2},\qquad
                 k_{\xv, \yv} = \sqrt{\frac{4 \xv_r \yv_r}{(\xv_r + \yv_r)^2 + (\xv_z - \yv_z)^2}}.
\end{align*}
The form $l$ contains the contribution from the coils, and is given by
\begin{align}
    l(I, v)&=
               {\sum_{i=1}^{N_I} \frac{I_i}{|\Omega_{{\rm C}_{i}}|} \int_{\Omega_{{\rm C}_{i}}} v \; dr dz. }
\end{align}
%
%

In order to solve~\eqref{eq:nonlin}, we introduce a Newton iteration for
$\varphi^{k+1} = \psi^{k+1} - \psi^k$. This iteration can be summarized as
the linearized system of equations for $\varphi^{k+1}$
\begin{align}
  &d_\psi(\varphi^{k+1}) [{\rm a}(\psi^k, v_i)]
    =
  l(I, v_i)
    - {\rm a}(\psi^k, v_i), \qquad i=1,\dots, N_{\rm dof}, \label{eq:solveme}
\end{align}
where $d_{\psi}(\varphi)$ represents the Gateaux semiderivative in the direction of $\varphi$ where,
for a smooth function $g$, is defined by
\begin{align}
  d_{\psi}(\varphi) [g(\psi)] = \lim_{\epsilon \rightarrow 0}\frac{g(\psi+\epsilon \varphi) - g(\psi)}{\epsilon}
  = \frac{d}{d\epsilon}g(\psi+\epsilon \varphi) \bigg{|}_{\epsilon=0}
\end{align}
The Gateaux semiderivative of the operator $a$ is given by
\begin{align}
  d_\psi(\varphi) [a(\psi, v)] =
  d_\psi(\varphi) [a_{\Delta^*}(\psi, v)]
  + d_\psi(\varphi) [a_{\rm plasma}(\psi, v)]
  + d_\psi(\varphi) [a_{\rm far-field}(\psi, v)].
\end{align}
Due to linearity, $d_\psi(\varphi) [a_{\Delta^*}(\psi, v)] = a_{\Delta^*}(\varphi, v)$
and $d_\psi(\varphi) [a_{\rm far-field}(\psi, v)] = a_{\rm far-field}(\varphi, v)$.
The Gateaux semiderivative of $a_{\rm plasma}$ is given by
\begin{align}
  d_\psi(\varphi) [a_{\rm plasma}(\psi, v)]
  =\;&  - \int_{\Omega_{\rm p}(\psi)} d_{\psi}(\varphi)\left[ \left( r p'(\psi) + \frac{1}{\mu r} f(\psi) f'(\psi) \right) \right] v \; dr dz \nonumber\\
     &- \int_{\partial \Omega_{\rm p}(\psi)} \bigg{(} r p'(\psi) + \frac{1}{\mu r} f(\psi) f'(\psi) \bigg{)}
       |\grad \psi |^{-1} (\varphi(r_{\rm x}(\psi), z_{\rm x}(\psi)) - \varphi) \, v \; d s.
       \label{eq:plasmaterm}
\end{align}
See Section~\ref{sec:aplasma} for the derivation.

The derivatives $d_\psi(\varphi)[p'(\psi)]$ and $d_\psi(\varphi)[f(\psi)f'(\psi)]$ are nontrivial and depend on the choice for
representing these functions.
We compute these derivatives analytically using shape calculus, as introduced in~\cite{heumann2015quasi} for this problem.


\paragraph{Luxon and Brown Model}

The semiderivatives of the functions defined in~\eqref{eq:lb} are
\bse
\begin{align}
  d_\psi(\varphi)[p'(\psi)] &= \alpha \frac{\beta}{r_0} \gamma \delta (1 - \psi_N(\psi)^{\delta})^{\gamma-1} \psi_N(\psi)^{\delta-1} d_\psi(\varphi)[\psi_N(\psi)]
  \\
  d_\psi(\varphi)[f(\psi)f'(\psi)]
  &= \alpha (1 - \beta) \mu r_0 \gamma \delta (1 - \psi_N(\psi)^{\delta})^{\gamma-1} \psi_N(\psi)^{\delta-1} d_\psi(\varphi)[\psi_N(\psi)]
\end{align}
\ese
where
\begin{align}
  d_\psi(\varphi)[\psi_N(\psi)] = \frac{\varphi-   (1 - \psi_N(\psi)) \varphi(r_{\rm ma}(\psi), z_{\rm ma}(\psi)) - \psi_N(\psi)\varphi(r_{\rm x}(\psi), z_{\rm x}(\psi))}{\psi_{\rm x} - \psi_{\rm ma}}.
  \label{eq:shapecalculus}
\end{align}
The nontrivial derivation of~\eqref{eq:shapecalculus} involves shape calculus and is provided in Section~\ref{sec:semiempirical}.

\smallskip

\paragraph{Taylor Equilibrium}

The right-hand-side terms are given by
\begin{align}
  f(\psi)f'(\psi) = \alpha \left[ f_{\rm x} + \alpha (\psi - \psi_{\rm x})\right].
\end{align}
The relevant semiderivative is given by
\begin{align}
  d_\psi(\varphi)[f(\psi)f'(\psi)]
  &= \alpha^2 (1 - d_\psi(\varphi) \psi_{\rm x}) = \alpha^2 (1 - \varphi(r_{\rm x}(\psi), z_{\rm x}(\psi))).
\end{align}
Note that the result $d_\psi(\varphi) \psi_{\rm x} = \varphi(r_{\rm x}(\psi), z_{\rm x}(\psi))$, which is derived in Section~\ref{sec:semiempirical}, is used to simplify the result.

\smallskip


\paragraph{15MA ITER baseline case}

We have
\begin{align}
  f(\psi) f'(\psi) = \frac{1}{\psi_{\rm ma} - \psi_{\rm x}} \alpha (f_{\rm x} + \alpha S_f(\psi_N(\psi))) S_f'(\psi_N(\psi)).
\end{align}

Using the results from Section~\ref{sec:semiempirical}, the semiderivative for $ff'$ is given by
\begin{align}
  d_\psi(\varphi) [f(\psi)f'(\psi)]
  =\;& \frac{1}{\psi_{\rm ma} - \psi_{\rm x}} \alpha^2 S_f'(\psi_N(\psi)))^2 d_\psi(\varphi)[\psi_N(\psi)] \\
     &+ \frac{1}{\psi_{\rm ma} - \psi_{\rm x}} \alpha (f_{\rm x} + \alpha S_f(\psi_N(\psi))) S_f''(\psi_N(\psi)) d_\psi(\varphi)[\psi_N(\psi)] \nonumber \\
     &-\frac{\varphi_{\rm ma} - \varphi_{\rm x}}{(\psi_{\rm ma} - \psi_{\rm x})^2} \alpha (f_{\rm x} + \alpha S_f(\psi_N(\psi))) S_f'(\psi_N(\psi)).\nonumber
\end{align}
The semiderivative for $p'$ is given by
\begin{align}
  d_\psi(\varphi) [p'(\psi)] &= S_{p'}'(\psi_N(\psi)) d_\psi(\varphi)[\psi_N(\psi)].
\end{align}

\smallskip

Due to the presence of $\varphi_{\rm x}$ and $\varphi_{\rm ma}$ in \eqref{eq:plasmaterm},
the discrete matrix formed by building $d_\psi(\varphi) [a_{\rm plasma}(\psi, v)]$
contains nonzeros in columns corresponding to the indices of the magnetic axis and x-point
along rows corresponding to indices in $\Omega_p(\psi)$. We note that the resulting matrix is not symmetric.

\subsection{Location of $\xv_{\rm ma}$ and $\xv_{\rm x}$ and Determination of $\Omega_{\rm p}$}
\label{sec:separatrix}

In this section, the algorithms for locating the magnetic axis ($\xv_{\rm ma}$) and magnetic x-point ($\xv_{\rm x}$)
and determining the elements that contain $\Omega_{\rm p}$ are described.
We assume that the connectivity of adjacent mesh vertices
is available in a convenient data structure, such as a map or dictionary.
Additionally, we assume that each vertex in $\Omega_{\rm L}$ or on
$\partial \Omega_{\rm L}$ is marked and collected into a list.
First, all of the marked vertices in $\Omega_{\rm L}$ and $\partial \Omega_{\rm L}$
are searched to determine the maximal and minimal values for $\psi$
in addition to any saddle points in $\psi$.
$\xv_{\rm ma}$ is set to be the location of the vertex corresponding to the minimum $\psi$ value.
If one or more vertices are found to be saddle points, $\xv_{\rm x}$ is set to be the location
of the saddle point vertex with a $\psi$ value closest to $\psi_{\rm ma}$.
Otherwise, $\xv_{\rm x}$ is set to be the location of the vertex corresponding to the maximum $\psi$ value.

Saddle points are determined using the following algorithm.
Let $\psi_o$ be the function value at a candidate vertex
and let $\psi_{1}, \dots, \psi_{N_v}$ be function values of adjacent vertices ordered clock-wise.
The candidate vertex is marked as a saddle point if the sequence $\psi_1-\psi_o, \dots, \psi_{N_v}-\psi_o, \psi_1-\psi_o$
contains 4 or more sign changes.

After determining $\xv_{\rm ma}$ and $\xv_{\rm x}$, a tree search is performed to determine the vertices
belonging to or adjacent to $\Omega_{\rm p}$. A dynamic queue is initially populated with the vertex at $\xv_{\rm ma}$.
Until the queue is empty, vertices in the queue are removed from the queue and its adjacent vertices are
investigated. If a neighboring vertex has a $\psi$ value in between $\psi_{\rm ma}$ and $\psi_{\rm x}$,
then it is marked to be inside of $\Omega_{\rm p}$ and added to the queue.
Otherwise, it is marked to be adjacent to $\Omega_{\rm p}$.
Domain integrals over $\Omega_{\rm p}$ are evaluated over each element containing vertices in or
adjacent to $\Omega_{\rm p}$. When the vertex is adjacent to $\Omega_{\rm p}$, the integrand is set to be zero.

\subsection{Adaptive Mesh Refinement}
\label{sec:amrdesc}
Adaptive mesh refinement (AMR) is applied to increase the mesh resolution in areas of high priority
according to the Zienkiewicz-Zhu (ZZ) error estimation procedure~\cite{ZZ1,ZZ2}.
An element is marked for refinement when its local error estimate is greater than
a globally defined fraction of the total error estimate.
More specifically, we divide the mesh into two regions, $\Omega_{\rm L}$, which contains the
plasma region, and $\Omega \backslash \Omega_{\rm L}$, which have their own respective refinement thresholds.

We explore two approaches. In the first approach,  we use the computed solution $\psi$ and the elliptic operator from the Grad--Shafranov equation in the error estimate.
Specifically, the error of the ZZ estimator at an element $K$ can be estimated as
\begin{align*}
\eta_K^2(\psi_h) = \int_K\|G[\psi_h]-\tilde\nabla \psi_h\|^2 \; drdz
\end{align*}
where $G[\psi_h]$ denotes the smoother version of the flux, and $\tilde\nabla \psi_h$ is a numerical interpolated flux.
Here $G[\psi_h]$ is solved from the diffusion operator $\frac{1}{\mu r} \Delta^*$,
and $\tilde\nabla:=\frac{1}{\mu r}[\partial_r, \, \partial_z]^T$.
This is a well defined error estimator for the proposed problem, since the problem is of an elliptic nature.
However, this error estimator may prioritize the solution of $\psi$ instead of $\tilde\nabla\psi$, latter of which is more interesting as it becomes
the components of $B$.

To accommodate the need of resolving more interesting physics, we also consider using the toroidal $B$ field as the input in the ZZ estimator.
We first reconstruct the toroidal $B$ field from $f(\psi)$ and then use this as the input in the ZZ estimator.
This will naturally prioritize the jump in the toroidal current of the equilibrium.
Note that this approach leads to a feature-based indictor, which in practice can work well for the need.

\section{Optimization}
\label{sec:opt}
Consider a curve $\Gamma_{\rm p}$ describing the shape of a desired plasma boundary.
The goal is to find the poloidal magnetic flux $\psi(r, z)$
coil currents, $I_j$, $j=1,\dots, M$,
and scaling parameter, $\alpha$,
such that some objective function penalizing the distance between
$\partial \Omega_{\rm p}$ and $\Gamma_{\rm p}$,
\begin{align}
  g(\psi) = d(\partial \Omega_{\rm p}(\psi), \Gamma_{\rm p}),
\end{align}
is minimized subject to the following constraints.
The first constraint enforces that $\psi$ is a solution to the Grad--Shafranov equations,
which are summarized symbolically by the system of equations,
\begin{align}
  B(\psi) - F(I) = 0.
\end{align}
In addition to the above constraint, we enforce the following constraint to set the plasma current equal to
a desired value of $I_{\rm p}$,
\begin{align}
  {\int_{\Omega_{\rm p}(\psi)}   \left(r p'(\psi) + \frac{1}{\mu r} {ff'}(\psi)\right)   \; dr dz} = I_{\rm p}.
\end{align}
%
For well-posedness, we add the Tikhonov regularization term,
\begin{align}
  h(I) = \frac{1}{2}\sum_{j=1}^{N_I} w_j I_j^2,
\end{align}
which penalizes solutions with large currents.
There are $N_I$ independent currents that can be controlled.
The optimization problem can be written as
\begin{align*}
  \min_{\psi, \, \alpha, \,  I} \qquad&g(\psi) + h(I), \\
  \text{s.t.} \qquad& B(\psi) - F(I) = 0, \\
                    &  {\int_{\Omega_{\rm p}(\psi)}   \left(r p'(\psi) + \frac{1}{\mu r} {ff'}(\psi)\right)   \; dr dz} = I_{\rm p}.
\end{align*}
We now discuss the choice of objective function.
Points $x_1, \dots, x_{N_c}$ are chosen from a reference plasma boundary curve.
Since $\psi_N=1$ defines the contour that
contains the x-point, we define our objective to penalize when the solution along the reference curve is not
on the contour of the x-point,
\begin{align}
  g(\psi) = \frac{1}{2} \sum_{i=1}^{N_c} (\psi_N(\psi(x_i)) - 1)^2.
\end{align}
\das{Note that this objective is different from that in~\cite{heumann2015quasi},
where the objective was defined as the quadratic functional, equivalent to 
$g_{\rm alt}(\psi) = \frac{1}{2} \sum_{i=2}^{N_c} (\psi(x_i) - \psi(x_1))^2$.
The important distinction between $g_{\rm alt}(\psi)$ and $g(\psi)$ is that $g_{\rm alt}(\psi)$ promotes the field $\psi$
to have an isoline on the control points
while $g(\psi)$ promotes the field $\psi$ to have an isoline \textit{and the} separatrix on the control points.
Note that the separatrix is computed using the algorithm described in Section~\ref{sec:separatrix}.
}

\subsection{Fully Discretized System}

\das{The discussion of Section~\ref{sec:opt} so far focuses on the continuous version of the problem.
  The following will describe the discretized version of the continuous problem
  in terms of finite element coefficients.}
Consider the decomposition of $\psi$ into linear algebraic degrees of freedom, $y$.
Let $y_n$ be the coefficient for the
$n^{\rm th}$ shape basis function $v_{n}(r, z)$. Then
\begin{align}
  \psi(r, z) = \sum_{n=1}^{N_{\rm dof}} y_n  v_n(r, z).
\end{align}
We redefine the current variables as $u_j = I_j$ for $j=1,\dots, N_I$.
In terms of linear algebraic variables, we rewrite the objective function as
\begin{align*}
  \min_{\psi, I} \qquad&G(y) + R(u), \\
  \text{s.t.} \qquad& B(y, \alpha) - F u = 0, \\
                       & C(y, \alpha) - I_{\rm p} = 0,
\end{align*}
where
$G(y), R(u)\in\mathbb{R}$ represent the objective function and regularization,
$B(y, \alpha) \in\mathbb{R}^{N_{\rm dof}}$ represents the solution operator,
$F\in\mathbb{R}^{N_{\rm dof}\times N_I}$ represents the coil contributions,
and $C(y, \alpha)\in\mathbb{R}$ represents the plasma current.
The regularizer, $R(u)$, is defined as
\begin{align}
  R(u) = \frac{1}{2} u^T H u.
  \label{eq:reg}
\end{align}
A simple choice is to choose $H = \epsilon I$, where $\epsilon$ is small (e.g. $\epsilon\sim 10^{-12}$).

The Lagrangian for the problem is given by
\begin{align}
  \mathcal{L} = G(y) + R(u) + p^T (B(y, \alpha) - F u) + \lambda ( C(y, \alpha) - I_{\rm p})
\end{align}
where $p$ and $\lambda$ are Lagrange multipliers.
The fixed point $(y^*, u^*, \alpha^*, p^*, \lambda^*)$
is the solution to the equations
\bse
\label{eq:nonlinequations}
\begin{align}
  G_y(y) + B_y(y, \alpha)^T p + C_y(y, \alpha) \lambda &= 0, \label{eq:n1} \\
  Hu - F^T p &= 0, \\
  B(y, \alpha) - F u &= 0, \label{eq:n2}\\
  B_\alpha(y, \alpha)^T p + C_{\alpha}(y, \alpha) \lambda &= 0, \\
  C(y, \alpha) &= I_{\rm p}.
\end{align}
\ese

We introduce a Newton iteration to address~\eqref{eq:nonlinequations}.
Let $(y^n, u^n, p^n, \alpha^n, \lambda^n)$ represent the solution state at iteration $n$
and let $\Delta$ represent a forward difference operator (i.e., $\Delta y^n := y^{n+1} - y^n$).
\das{The Newton method compute $\Delta y^n$ by solving the following linearized system,}
\begin{align}
  &\left[\begin{array}{ccccc}
          G_{yy}(y^n)&0&B_y(y^n, , \alpha^n)^T&0&C_y(y, \alpha) \\
           0&H&-F^T&0&0 \\
           B_y(y^n, \alpha^n) & -F & 0 & B_\alpha(y, a)&0 \\
           0 & 0 & B_\alpha(y, \alpha)^T & 0&C_\alpha(y, \alpha)  \\
          C_y(y, \alpha)^T & 0 & 0 & C_\alpha(y, \alpha) & 0
        \end{array}\right]
  \left[\begin{array}{c}
          \Delta y^{n} \\
          \Delta u^{n} \\
          \Delta p^{n} \\
          \Delta \alpha^{n} \\
          \Delta \lambda^{n}
        \end{array}\right]=
    \left[\begin{array}{c}
            b_1^n\\b_2^n\\b_3^n\\b_4^n\\b_5^n
          \end{array}\right],\label{eq:5x5system}
\end{align}
where
\begin{align}
  \left[\begin{array}{c}
          b_1^n\\b_2^n\\b_3^n\\b_4^n\\b_5^n
        \end{array}\right]=
    -\left[\begin{array}{c}
           G_y(y^n) + B_y(y^n, \alpha^n)^T p^n + C_y(y^n, \alpha^n) \lambda^n \\
           H u^n - F^T p^n \\
           B(y^n, \alpha^n) - Fu^n \\
           B_\alpha(y^n, \alpha^n)^T p^n + C_{\alpha}(y^n, \alpha^n) \lambda^n \\
           C(y^n, \alpha^n) - I_{\rm p}
                                                           \end{array}\right].
\end{align}
We eliminate $\Delta u^n$, $\Delta \alpha^n$, and $\Delta \lambda^n$ algebraically,
\das{
\begin{align}
  \Delta u^n = H^{-1} b_2 + H^{-1} F^T \Delta p^n, \qquad
  \Delta \alpha^n = \frac{b_5 - C_y^T \Delta y^n}{C_\alpha}, \qquad
  \Delta \lambda^n = \frac{b_4 - B_{\alpha}^T \Delta p^n}{C_\alpha},
  \label{eq:algebra}
\end{align}
}
to arrive
at a block system given by
%
    %
    %
%
\begin{align}
  \left[\begin{array}{cc}
          B_y - \frac{1}{C_\alpha}B_{\alpha} C_y^T & - F H^{-1} F^T \\
          G_{yy} & B_y^T - \frac{1}{C_\alpha} C_y B_{\alpha}^T
        \end{array}\right]
  \left[\begin{array}{c}
          \Delta y \\ \Delta p
        \end{array}\right]
  =
  \left[\begin{array}{c}
          b_3 + F H^{-1} b_2 - \frac{1}{C_\alpha} B_\alpha b_5  \\
          b_1 - \frac{1}{C_\alpha} C_y b_4
        \end{array}\right].
  \label{eq:2x2system}
\end{align}
For notational convenience, the functional dependence of the operators and the superscripts denoting the iteration
number are omitted.
\das{Therefore, our strategy involves solving~\eqref{eq:2x2system} for $\Delta y, \Delta p$ using
  advanced preconditioning techniques
  and using~\eqref{eq:algebra} to obtain $\Delta u, \Delta \alpha, \Delta \lambda$.}
We note that the original \das{(equivalent)} system~\eqref{eq:5x5system} is symmetric while the system~\eqref{eq:2x2system} after the block factorization is non-symmetric.
Although the system  \eqref{eq:2x2system} can be easily symmetrized by reordering $\Delta y$ and $\Delta p$, this does not help the preconditioning strategy that will be discussed in the next section.

\subsection{Computation of Objective Function and its Derivatives}

\das{The solution of}~\eqref{eq:2x2system} requires the computation of gradient and Hessian of the objective function.
In terms of the linear algebraic variables, the objective function is given by
\begin{align}
  G(y) = \frac{1}{2}\sum_{i=1}^{N_c} (\bar{y}_i - 1)^2
\end{align}
where $\bar{y}_i$ represents the discrete analog to $\psi_N(\psi(x_i))$.
For $H^1(\Omega)$ functions represented on 2D triangular elements, off-nodal function values can be expressed as an
interpolation of the nodal values.
For each control point, $x_i$, the bounding vertices are located and their indices are labeled as $J_1^{(i)}, J_2^{(i)}, J_3^{(i)}$.
Additionally, corresponding interpolation coefficients, $\alpha_{J_1^{(i)}}, \alpha_{J_2^{(i)}}, \alpha_{J_3^{(i)}},$ are computed.
We therefore define $\bar{y}_i$ as
\begin{align}
  \bar{y}_i = \frac{\sum_{m=1}^{3} \alpha_{J_m^{(i)}} y_{J_m^{(i)}} - y_{\rm ma}}{y_{\rm x} - y_{\rm ma}},
\end{align}
where $y_{\rm ma}$ and $y_{\rm x}$ correspond to the nodal values of the magnetic axis and x-point.
The gradient and Hessian of the objective function are defined through
\bse
\begin{align}
  \frac{\partial }{\partial y_m} G(y) &= \sum_{i=1}^{N_c} (\bar{y}_i - 1) \frac{\partial \bar{y}_i}{\partial y_m}, \\
  \frac{\partial^2 }{\partial y_m \partial y_n} G(y) &=
  \sum_{i=1}^{N_c} \left( (\bar{y}_i - 1) \frac{\partial^2 \bar{y}_i}{\partial y_m \partial y_n}
  + \frac{\partial \bar{y}_i}{\partial y_m} \frac{\partial \bar{y}_i}{\partial y_n} \right).
\end{align}
\ese
The above require the definitions of the following derivatives.
\bse
\begin{align}
  \frac{\partial \bar{y}_i}{\partial y_m} &=
  \frac{1}{y_{\rm x} - y_{\rm ma}}
  \left(\sum_{m=1}^{3} \alpha_{J_m^{(i)}} \frac{\partial y_{J_m^{(i)}}}{\partial y_m} - \frac{\partial y_{\rm ma}}{\partial y_m}\right)
  - \frac{1}{y_{\rm x} - y_{\rm ma}} \left(\frac{\partial y_{\rm x}}{\partial y_{m}} - \frac{\partial y_{\rm ma}}{\partial y_{m}}\right) \bar{y}_i, \\
  \frac{\partial^2 \bar{y}_i}{\partial y_m \partial y_n} &=
  \frac{1}{y_{\rm ma} - y_{\rm x}} \left(  \left(\frac{\partial y_{\rm x}}{\partial y_{n}} - \frac{\partial y_{\rm ma}}{\partial y_{n}}\right) \frac{\partial \bar{y}_i}{\partial y_m}
  + \left(\frac{\partial y_{\rm x}}{\partial y_{m}} - \frac{\partial y_{\rm ma}}{\partial y_{m}}\right) \frac{\partial \bar{y}_i}{\partial y_n}\right)
\end{align}
\ese

\section{Preconditioning Approach}
\label{sec:prec}

At each iteration, the system is given by
\begin{align}
  \left[\begin{array}{cc}
          B & A \\ C & B^T
        \end{array}\right]
  \left[\begin{array}{c}
          \Delta p \\ \Delta y
        \end{array}\right]
  =
  \left[\begin{array}{c}
          c_1 \\ c_2
        \end{array}\right],
  \label{eq:2x2sys}
\end{align}
where the block operators are given by
\begin{align}
  A = G_{yy}, \qquad
  B = B_y^T - \frac{1}{C_\alpha} C_y B_{\alpha}^T, \qquad
  C = - F H^{-1} F^T,
\end{align}
and the right hand side block vectors are defined to be
\begin{align}
  c_1 = b_1 - \frac{1}{C_\alpha} C_y b_4,\qquad
  c_2 = b_3 + F H^{-1} b_2 - \frac{1}{C_\alpha} B_\alpha b_5.
\end{align}
$A$ and $-C$ are both symmetric positive semi-definite
and $B$ and $B^T$ are both invertible.
The operator $B$ incorporates the contribution from the elliptic operator, plasma terms, and boundary conditions
in $B_y^T$
and a rank-one perturbation in $-\frac{1}{C_\alpha} C_y B_{\alpha}^T$ that is dense over the plasma degrees of freedom.
The operator $C$ is dense over the coil degrees of freedom and the
operator $A$ only contains non-zero entries nearby the plasma control points.

At each Newton iteration,~\eqref{eq:2x2sys} is solved using FGMRES~\cite{fgmres93}.
The following discussion summarizes the preconditioning approach.
In developing a preconditioner, our hypothesis is that the problem is dominated primarily by the elliptic part of
the diagonal operators. This motivates a block diagonal preconditioner using algebraic multigrid (AMG).
\begin{align}
  \text{Block Diagonal PC:} \qquad \mathcal{P}^{-1}_{\rm BD} &=
  \left[\begin{array}{cc}
          {\rm AMG}(B_y^T) & 0 \\ 0 & {\rm AMG}(B_y)
        \end{array}\right].
\label{eq:pca}
\end{align}
\das{Here, AMG represents a sub-matrix preconditioner that applies algebraic multigrid to the input operator
  using a specified cycle type (V or W-cycle) and number of iterations.}
This \das{block} preconditioner uses \das{two calls of the AMG preconditioner on the sub-matrices} per Krylov solver iteration,
\das{one for $B_y^T$ and one for $B_y$}.
Another level of approximation involves incorporating
one of the off-diagonal blocks.
\begin{align}
  \text{Block Upper Triangular PC:} \qquad \mathcal{P}^{-1}_{\rm BUT} &=
                     \left[\begin{array}{cc}
                             {\rm AMG} \left(B_y^T\right) & D \\ 0 & {\rm AMG} \left(B_y\right)
                           \end{array}\right],\label{eq:pcb}\\
  \text{Block Lower Triangular PC:} \qquad \mathcal{P}^{-1}_{\rm BLT} &=
                     \left[\begin{array}{cc}
                             {\rm AMG} \left(B_y^T\right) & 0 \\ E & {\rm AMG} \left(B_y\right)
                           \end{array}\right].            \label{eq:pcc}
\end{align}
Here, $D = -{\rm AMG}(B_y^T)\, A\, {\rm AMG}(B_y)$
and $E= -{\rm AMG}(B_y)\, C\, {\rm AMG}(B_y^T)$.
These preconditioners also require two \das{AMG calls of  ${\rm AMG}(B_y^T)$ and ${\rm AMG}(B_y)$} per Krylov solver iteration.
For each preconditioner, we consider their performance for the choice of AMG cycle type
(either V-cycle or W-cycle) and number of AMG iterations.

\begin{remark}
  The system can be solved by calling AMG as the preconditioner for the whole block system.
  \dasb{For example, this approach is employed to solve the 
    biharmonic equation in one of the examples from the hypre library~\cite{hypre}.
    The biharmonic equation is transformed into a system of two 
    coupled Laplace equations which shares similarities with~\eqref{eq:2x2sys}.
    However for our problem,
    calling AMG on the whole block system is found to be less efficient than the proposed block preconditioning strategy.}
\end{remark}


\paragraph{Inexact Newton Method}
An inexact Newton procedure based on~\cite{eisenstat96} is used to determine the relative tolerance
for the FGMRES solver. Let $e_n$ be the global error norm at step $n$ and let $\eta_{\max}$ be the maximum
relative tolerance, then we choose the relative tolerance at the next Newton iteration to be given by
\begin{align}
  \eta_n = \min\left(\gamma \left(\frac{e_n}{e_{n-1}} \right)^{\theta}, \eta_{\max} \right),
  \label{eq:inexact}
\end{align}
where $\theta = (1 + \sqrt{5}) / 2$ and $\gamma\approx 1$.

\section{Implementation} \label{sec:mfem}

Our algorithm for solving the Grad--Shafranov equation is implemented using the freely available
C++ library MFEM~\cite{mfem-web, anderson2021mfem}.
In this section, a few of the implementation details are discussed.
The geometry is meshed using triangular elements using the Gmsh software~\cite{gmsh}.
The linearized system of equations in~\eqref{eq:solveme} is implemented using a combination
of built-in and custom integrators.
In particular, the bilinear form associated with~\eqref{eq:termc} is built using a custom integrator
and the contribution to $B_y$ from~\eqref{eq:plasmaterm} is implemented by building a
custom sparse matrix.
The block system in~\eqref{eq:2x2sys}
and the preconditioners in~\eqref{eq:pca}--\eqref{eq:pcb}
are implemented as nested block operators.
Hypre's algebraic multigrid preconditioners~\cite{hypre}
are used to precondition the matrices $B_y$ and $B_y^T$
and MFEM's built in FGMRES implementation is used to solve the block system.
\das{The default values given by the MFEM's hypre wrapper are used except two \das{parameters}: the cycle type and 
  the number of multigrid iterations per preconditioner, which will be studied in the numerical section.}
\dasb{The smoothing operator and number of smoothing steps are chosen to be the default choices of the wrapper, which are l1-Gauss-Seidel and 1.}

MFEM's built in Zienkiewicz-Zhu error estimation procedure is used to refine the mesh at each AMR level.
The implementation of the solver in this paper and the scripts used to generate the results of this
paper are currently maintained in the \texttt{tds-gs} branch of MFEM at
\texttt{https://github.com/mfem/mfem/tree/tds-gs}.
During the preparation of this paper, we used MFEM version 4.5.3 and Hypre version 2.26.0.

\section{Results}
\label{sec:results}
In this section, we demonstrate our new solver using the three discussed approaches for handling the right-hand-side
of~\eqref{eq:govEq}.
First we perform a study to test the performance of the new preconditioner.
We then devise a case study to show the benefit of using adaptive mesh refinement.

\subsection{Preconditioner Performance}
\label{sec:precperf}
\begin{figure}[htb]
  \centering
  \raisebox{-.5\height}{\includegraphics[width=.66\textwidth, trim={3cm 0 3cm 0}, clip]{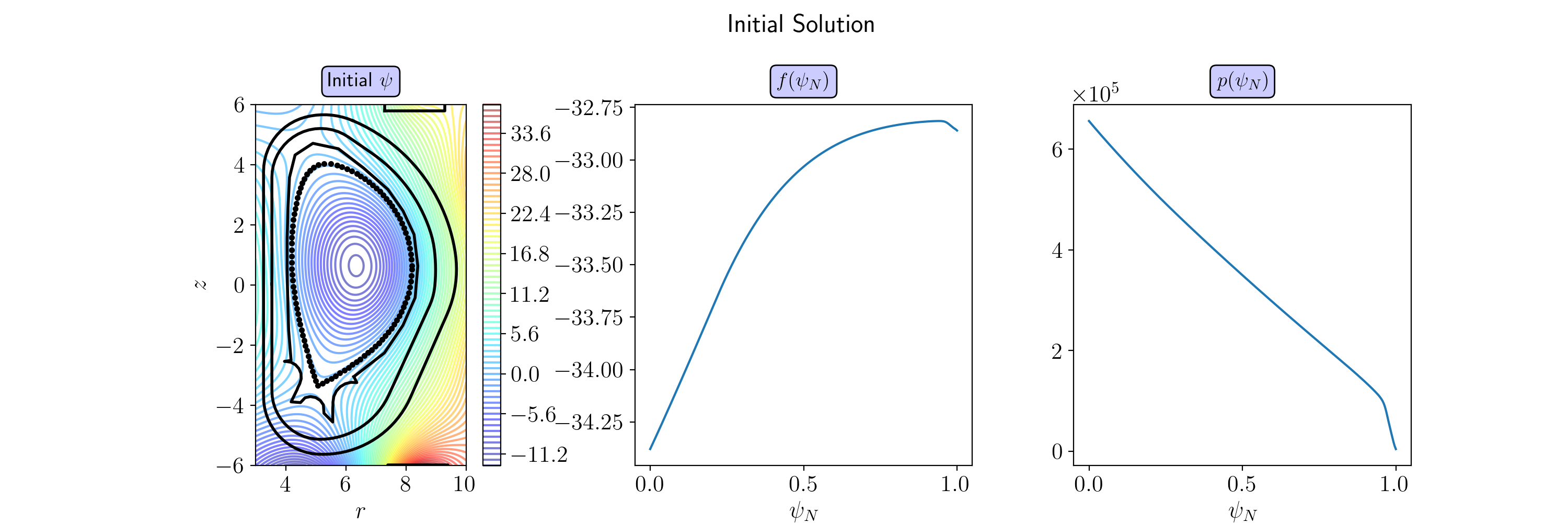}}
  \resizebox{.3\textwidth}{!}{
    \begin{tabular}{|c|c|c|c|}
    \hline
      \multicolumn{4}{|c|}{Initial Coil Currents} \\
    \hline
    \multicolumn{2}{|c|}{Center Solenoids} & \multicolumn{2}{c|}{Poloidal Flux Coils} \\
    \hline
    $I_1$ & 1.143284e+03 & $I_6$ & -4.552585e+06 \\
    \hline
    $I_2$ & -2.478694e+04 & $I_7$ & 3.180596e+06 \\
    \hline
    $I_3$ & -3.022037e+04 & $I_8$ & 5.678096e+06 \\
    \hline
    $I_4$ & -2.205664e+04 & $I_9$ & 3.825538e+06 \\
    \hline
    $I_5$ & -2.848113e+03 & $I_{10}$ & 1.066498e+07 \\
    \hline
          && $I_{11}$ & -2.094771e+07 \\
    \hline
  \end{tabular}
  }
  \caption{Details of the reference solution used for the initial guess.
    Left: initial field for $\psi$.
    Middle left: line plot of $f(\psi_N)$.
    Middle right: line plot of $p(\psi_N)$.
    Right: initial coil currents.
    \label{fig:initial}}
\end{figure}

We consider the Luxon and Brown model, Taylor state equilibrium, and the 15MA ITER baseline case.
For all cases, an initial guess was provided from a proposed ITER discharge at 15MA toroidal plasma current \cite{liu2015modelling},
which carries the ITER reference number ABT4ZL.
The initial solution for $\psi$ is provided in the domain $(r, z) \in [3, 10] \times [-6, 6]$
and shown in Figure~\ref{fig:initial}.
The initial guess for coil currents is shown in the table of Figure~\ref{fig:initial}.
100 control points along the plasma boundary are chosen for the objective function.
The magnetic permeability is chosen to be $\mu = 1.25663706144 \cdot 10^{-6}$ and the target
plasma current is set to $I_{\rm p}=15$ MA.
For all cases, the Newton \dasb{absolute} tolerance was set to $10^{-6}$.
\dasb{The residual is computed on a non-dimensionalized form of the nonlinear equations in
  \eqref{eq:2x2sys} that tames the large scale arising from $1/\mu$.
}
Problem specific solver settings are provided in Table~\ref{tab:settings}.
The regularization coefficient, $\epsilon$, maximum Krylov relative tolerance, $\eta_{\max}$, initial
uniform refinements, and total AMR levels were manually varied based on the problem.
\das{As $\epsilon$ affects the balancing between the objective function and regularization, it mainly
  affects the resulting solution and minimally impacts the preconditioner performance.
  $\eta_{\max}$ affects the number of FMGRES iterations to reach the tolerance and was chosen
  roughly to be the \dasb{largest} value that provided a \dasb{consistent} number of Newton steps.
  The initial number of uniform refinements affects the resolution of the initial mesh and the 
  number of AMR levels is the number of additional AMR refinements in the solver.
  We found that the Luxon and Brown and 15MA ITER profiles
  required a finer mesh in order to robustly reach convergence for the Newton solver.
  The number of AMR levels for each profile was chosen to balance total computational time
  of running all the cases to generate the results.
  }

\begin{table}[htb]
  \centering
  \caption{Solver settings chosen for each problem in Section~\ref{sec:precperf}.
    $\epsilon$ corresponds to the regularization coefficient in $H=\epsilon I$, defined in~\eqref{eq:reg}.
    $\eta_{\max}$ is the maximum Krylov relative tolerance used in~\eqref{eq:inexact}.
    Also reported is the number of uniform refinements performed on the initial mesh
    and the number of AMR levels used in the solver.
      \label{tab:settings}}
  \begin{tabular}{|c|c|c|c|}
    \hline
    & Luxon and Brown & Taylor State & 15MA ITER \\
    \hline
    $\epsilon$ & 1e-12 & 1e-12 & 1e-14 \\
    \hline
    $\eta_{\max}$ & 1e-6 & 1e-4 & 1e-6 \\
    \hline
    \multicolumn{1}{|c|}{initial uniform refinements} & 1 & 0 & 2 \\
    \hline
    \multicolumn{1}{|c|}{AMR levels} & 1 & 4 & 0 \\
    \hline
  \end{tabular}
\end{table}

Table~\ref{tab:iterations} summarizes the performance of the various preconditioner options chosen
in terms of number of Newton iterations and outer FGMRES iterations.
\das{In addition to choosing between the block diagonal, block upper triangular, and
  block lower triangular preconditioners, 
  the options that were varied were AMG iterations, which is the number of
  cycles performed on the matrix operator,
  and cycle type (V- or W-cycle).
  }
The number of Newton iterations on the initial mesh are reported in addition to the average number
of Newton iterations for each subsequent AMR level.
For the majority of cases considered, 4-8 Newton iterations are required to reach an \dasb{absolute} error tolerance of $10^{-6}$.
The notable exception includes applying multiple iterations of a W-cycle for the Luxon and Brown model
where performance degradation is observed for all preconditioners.
Excluding the cases with performance degradation, monotonic improvements in the average number of iterations
are seen by either moving from the block-diagonal PC to one of the block triangular PCs,
going from a V-Cycle to a W-Cycle, and
increasing the number of AMG iterations.
While performance is improved with more iterations, there are diminishing returns due to the extra cost of
additional AMG iterations.
We do not observe a significant difference in performance between the upper block triangular and the
lower block triangular PCs.

\begin{table}[hbt]
  \centering
  \caption{Summary of Newton and FGMRES iterations for the block system in~\eqref{eq:2x2sys}
    using the block diagonal, upper block triangular, and lower block triangular preconditioners
    for various combinations of AMG iterations and cycle type.
    The results are reported in the format $(a_1, \bar{a})\, b$ or $(a_1)\, b$, where $a_1$ is
    the number of Newton iterations on the initial mesh,
    $\bar{a}$ is the average number of Newton iterations on subsequent AMR levels
    and $b$ is the average number of FGMRES iterations.
    \label{tab:iterations}}
  Luxon and Brown

  \vspace{.1cm}

  \resizebox{.33\textwidth}{!}{
  \begin{tabular}{|c|c|c|}
    \hline
    \multicolumn{3}{|c|}{Block Diagonal PC} \\
    \hline
    {\small AMG \das{Itrs}} & V-Cycle & W-Cycle \\
    \hline
1 & (5, 2.0) 91.0 & (5, 2.0) 48.4 \\
3 & (5, 2.0) 57.7 & (5, 2.0) 40.3 \\
5 & (5, 2.0) 47.1 & (15, 2.0) 44.5 \\
10 & (5, 2.0) 37.4 & -- \\
    \hline
  \end{tabular}
  }\resizebox{.33\textwidth}{!}{
  \begin{tabular}{|c|c|c|}
    \hline
    \multicolumn{3}{|c|}{Upper Block Triangular PC} \\
    \hline
    {\small AMG Itrs} & V-Cycle & W-Cycle \\
    \hline
1 & (5, 2.0) 81.4 & (5, 2.0) 38.3 \\
3 & (5, 2.0) 47.0 & (5, 2.0) 32.4 \\
5 & (5, 2.0) 36.1 & (17, 2.0) 29.0 \\
10 & (5, 2.0) 26.9 & -- \\
    \hline
  \end{tabular}
  }\resizebox{.33\textwidth}{!}{
  \begin{tabular}{|c|c|c|}
    \hline
    \multicolumn{3}{|c|}{Lower Block Triangular PC} \\
    \hline
    {\small AMG Itrs} & V-Cycle & W-Cycle \\
    \hline
1 & (5, 2.0) 81.3 & (5, 2.0) 37.7 \\
3 & (5, 2.0) 47.0 & (5, 2.0) 38.1 \\
5 & (5, 2.0) 36.0 & (20, 2.0) 29.8 \\
10 & (5, 2.0) 26.4 & -- \\
    \hline
  \end{tabular}
  }

  \medskip

  Taylor State Equilibrium

  \vspace{.1cm}

  \resizebox{.33\textwidth}{!}{
  \begin{tabular}{|c|c|c|}
    \hline
    \multicolumn{3}{|c|}{Block Diagonal PC} \\
    \hline
    {\small AMG Itrs} & V-Cycle & W-Cycle \\
    \hline
1 & (5, 2.0) 69.5 & (5, 2.0) 38.2 \\
3 & (5, 2.0) 46.4 & (4, 2.0) 30.4 \\
5 & (4, 2.0) 41.1 & (4, 2.0) 28.3 \\
10 & (4, 2.2) 35.2 & (4, 2.0) 26.9 \\
    \hline
  \end{tabular}
  }\resizebox{.33\textwidth}{!}{
  \begin{tabular}{|c|c|c|}
    \hline
    \multicolumn{3}{|c|}{Upper Block Triangular PC} \\
    \hline
    {\small AMG Itrs} & V-Cycle & W-Cycle \\
    \hline
1 & (8, 2.0) 53.8 & (5, 2.0) 27.0 \\
3 & (5, 2.2) 36.4 & (5, 2.0) 18.8 \\
5 & (5, 2.2) 29.1 & (4, 2.0) 16.7 \\
10 & (5, 2.2) 22.7 & (4, 2.0) 14.6 \\
    \hline
  \end{tabular}
  }\resizebox{.33\textwidth}{!}{
  \begin{tabular}{|c|c|c|}
    \hline
    \multicolumn{3}{|c|}{Lower Block Triangular PC} \\
    \hline
    {\small AMG Itrs} & V-Cycle & W-Cycle \\
    \hline
1 & (5, 2.0) 59.4 & (5, 2.0) 27.5 \\
3 & (5, 2.2) 36.9 & (4, 2.0) 19.6 \\
5 & (4, 2.0) 28.8 & (4, 2.0) 17.4 \\
10 & (4, 2.2) 23.5 & (4, 2.0) 15.9 \\
    \hline
  \end{tabular}
}

  \medskip

  15MA ITER Case

  \vspace{.1cm}

  \resizebox{.33\textwidth}{!}{
  \begin{tabular}{|c|c|c|}
    \hline
    \multicolumn{3}{|c|}{Block Diagonal PC} \\
    \hline
    {\small AMG Itrs} & V-Cycle & W-Cycle \\
    \hline
1 & (7) 176.1 & (7) 67.7 \\
3 & (7) 107.3 & (7) 46.4 \\
5 & (8) 85.6 & (7) 40.6 \\
10 & (7) 63.9 & (7) 35.7 \\
    \hline
  \end{tabular}
  }\resizebox{.33\textwidth}{!}{
  \begin{tabular}{|c|c|c|}
    \hline
    \multicolumn{3}{|c|}{Upper Block Triangular PC} \\
    \hline
    {\small AMG Itrs} & V-Cycle & W-Cycle \\
    \hline
1 & (7) 166.3 & (7) 56.4 \\
3 & (7) 96.1 & (7) 34.6 \\
5 & (7) 74.4 & (7) 28.7 \\
10 & (7) 52.3 & (7) 23.9 \\
    \hline
  \end{tabular}
  }\resizebox{.33\textwidth}{!}{
  \begin{tabular}{|c|c|c|}
    \hline
    \multicolumn{3}{|c|}{Lower Block Triangular PC} \\
    \hline
    {\small AMG Itrs} & V-Cycle & W-Cycle \\
    \hline
1 & (7) 167.0 & (7) 56.7 \\
3 & (7) 96.6 & (7) 35.3 \\
5 & (7) 74.7 & (7) 29.6 \\
10 & (7) 52.9 & (7) 24.7 \\
    \hline
  \end{tabular}
  }
\end{table}

\begin{figure}[htb]
  \centering
  \includegraphics[width=\textwidth, trim={.5cm .5cm .5cm 0}, clip]{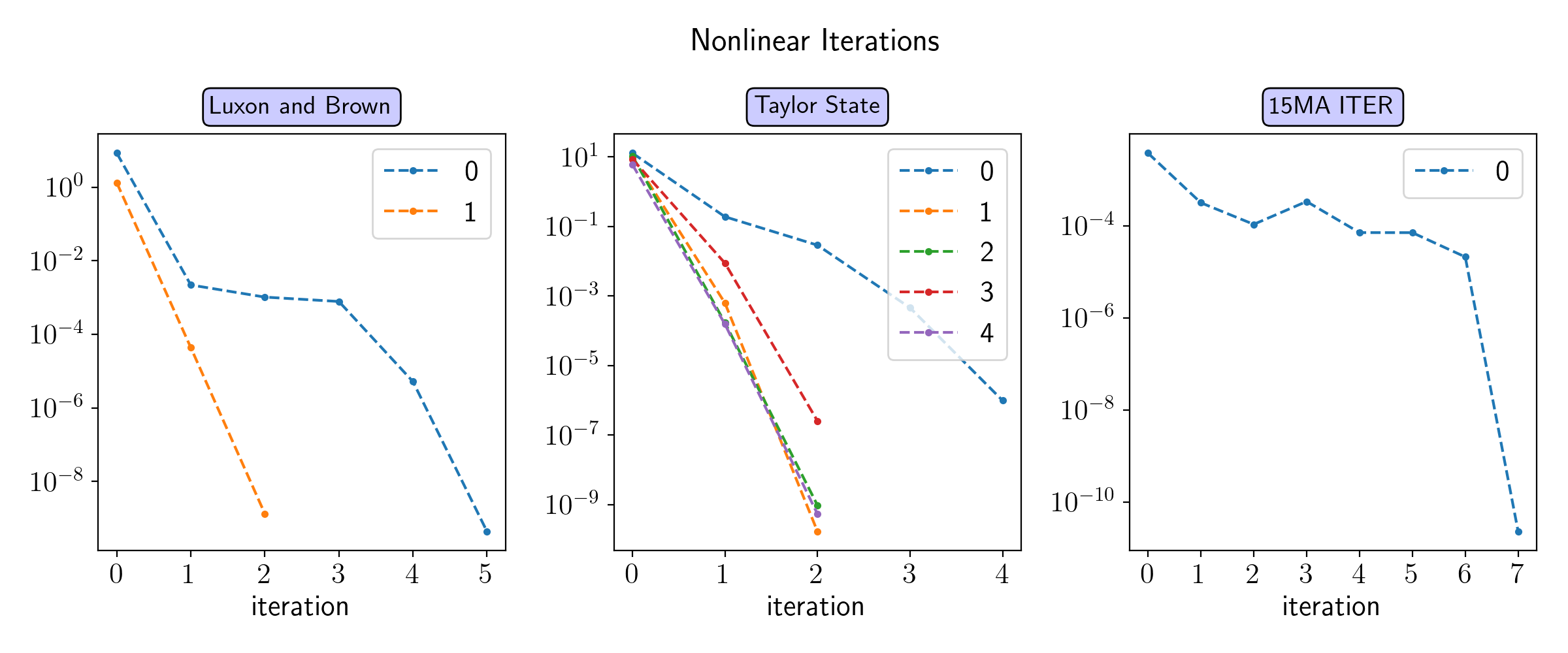}
  \caption{Newton iteration history \das{of the norm of the errors
    of equation~\eqref{eq:nonlinequations}}
    for each AMR refinement of  solution corresponding to the
    upper block triangular preconditioner with 10 AMG V-Cycle iterations for Luxon and Brown
    and 10 AMG W-Cycle iterations for the Taylor State and 15MA ITER.
    \label{fig:newtoniter}}
\end{figure}

\begin{figure}[htb]
  \centering
  \includegraphics[width=.8\textwidth, trim={0cm 0 0cm 0}, clip]{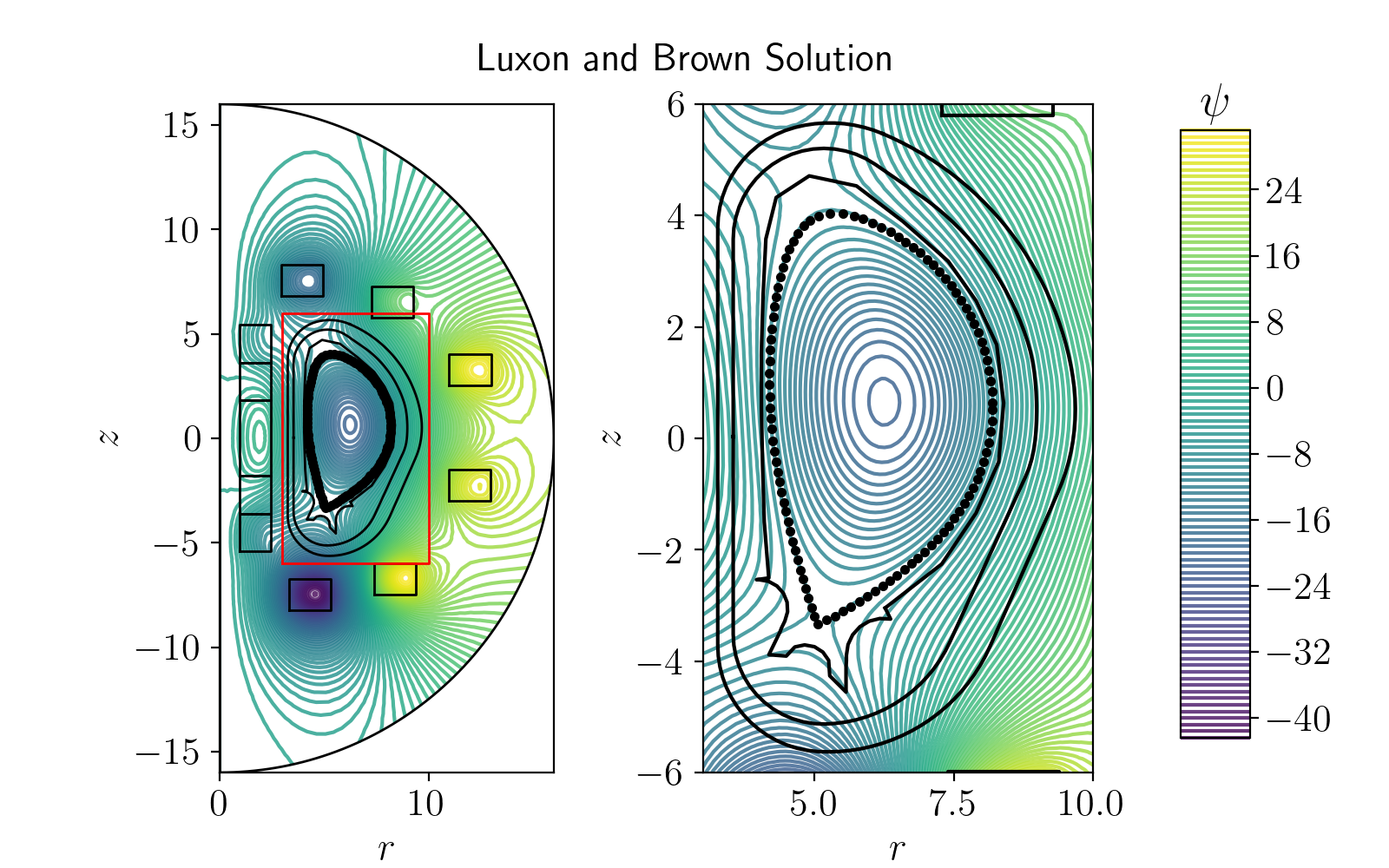}
    \caption{Luxon and Brown solution. Left: full domain. Right: zoomed in domain.
      \label{fig:lb}}
\end{figure}

\begin{figure}[htb]
  \centering
  \includegraphics[width=.8\textwidth, trim={0cm 0 0cm 0}, clip]{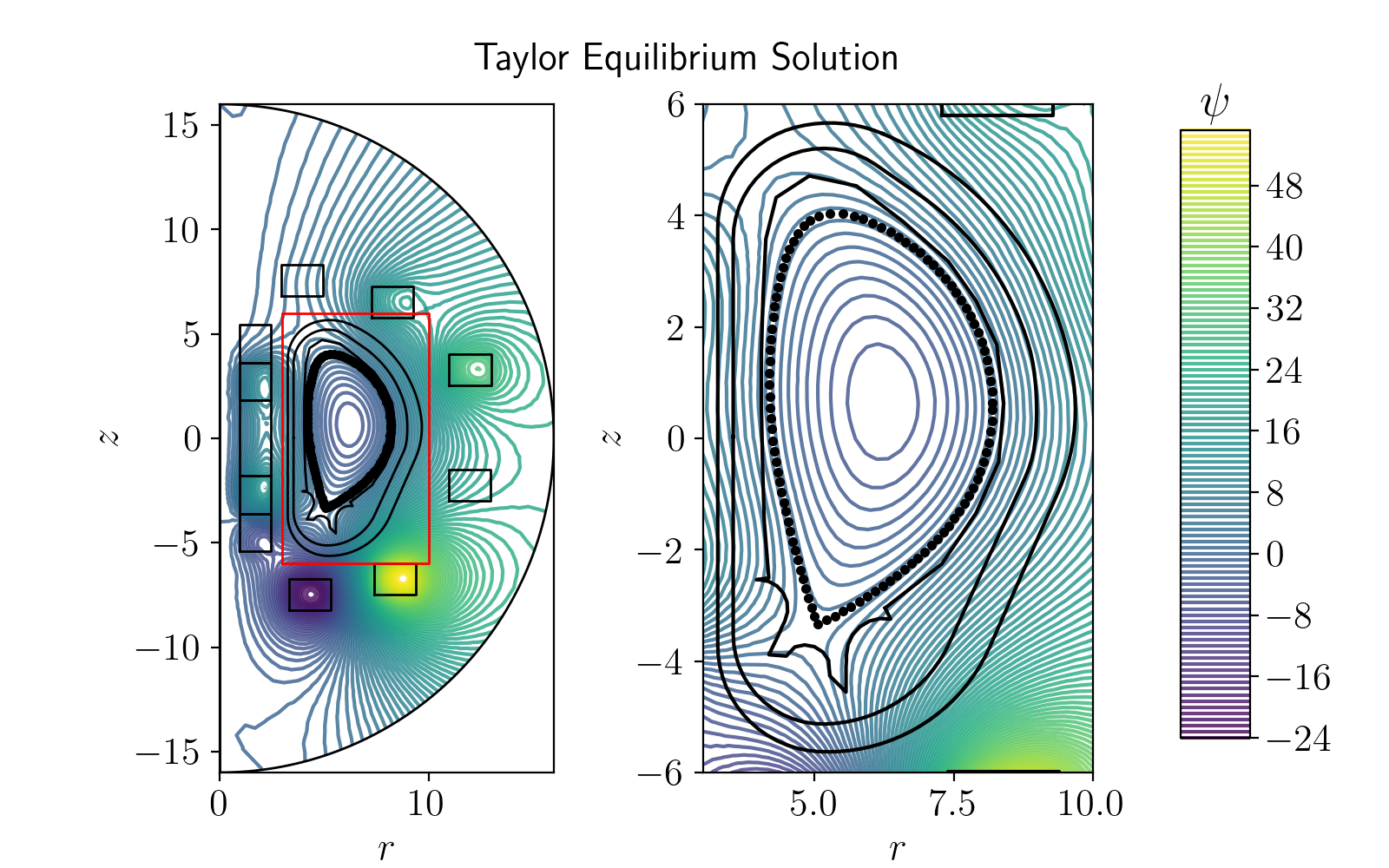}
    \caption{Taylor state equilibrium solution. Left: full domain. Right: zoomed in domain.
      \label{fig:taylor}}
\end{figure}

\begin{figure}[htb]
  \centering
  \includegraphics[width=.8\textwidth, trim={0cm 0 0cm 0}, clip]{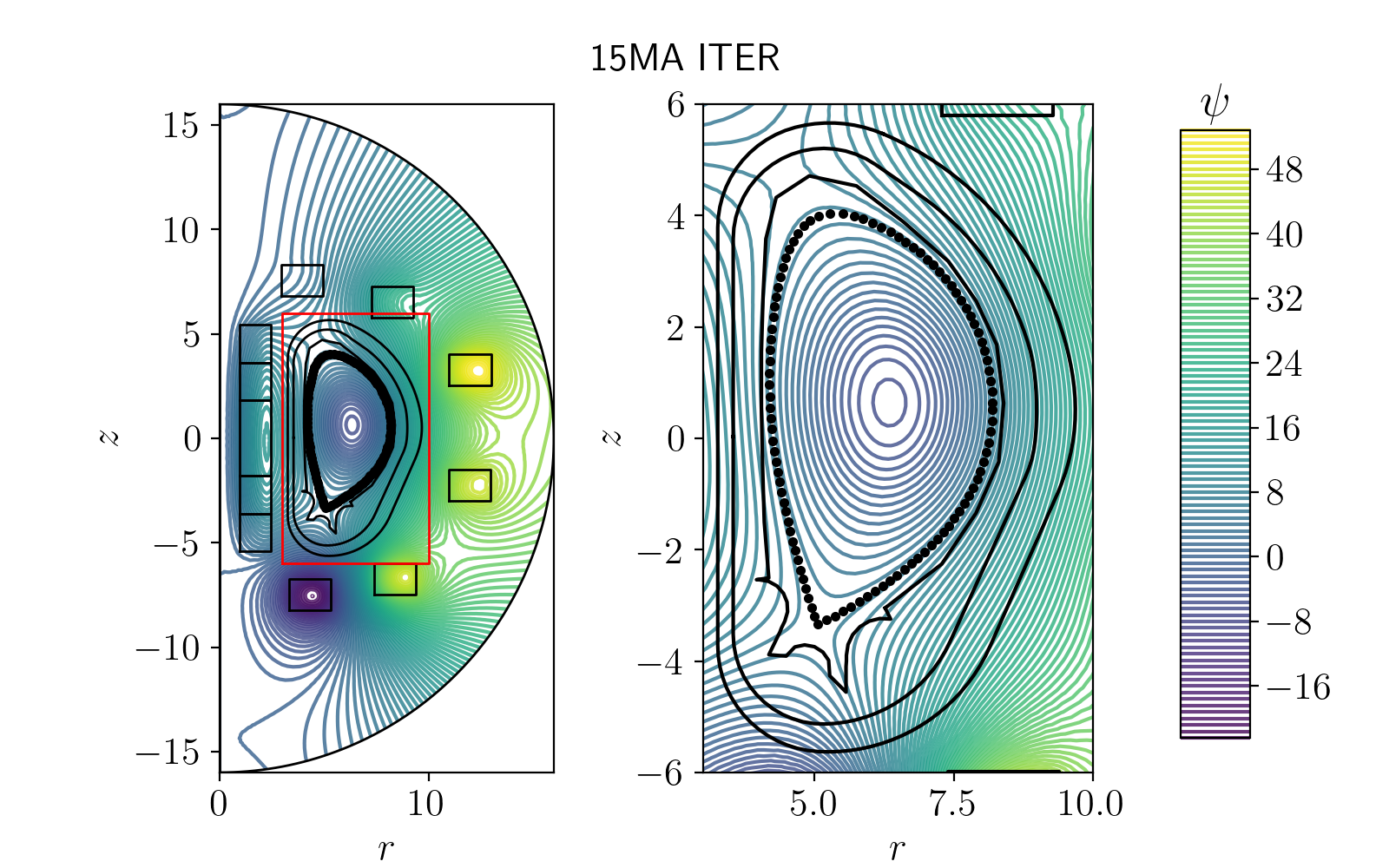}
    \caption{15MA ITER solution. Left: full domain. Right: zoomed in domain.
      \label{fig:empirical}}
\end{figure}

On subsequent AMR levels, an average of only 2 Newton iterations are generally required.
The convergence of the Newton iteration for a few selected cases is shown in Figure~\ref{fig:newtoniter}.
We note that the convergence rate of the latter AMR iteration becomes much better due to a better initial guess from the previous solve being used,
highlighting the power of AMR.
\dasb{Note that the Newton absolute tolerance is set to $10^{-6}$, therefore the iteration terminates when a smaller equation residual is achieved.
It is also evident from Figure~\ref{fig:newtoniter} that in all the studied cases, the relative change of residuals, defined as $\delta_r = \| \mathbf{r}\|/\|\mathbf{r}_0\|$, is $10^{-7}$ or much smaller.
Here we follow PETSc's definition of Newton convergence criteria {\tt SNESConvergedReason}, although other variations of Newton convergence criteria have been used in practice \cite{shadid2016scalable}.
When seeking a harder equilibrium, a combination of absolute and relative tolerance like the one used in \cite{shadid2016scalable} might be necessary to achieve converged solutions. 
}
We also note that the initial error for the 15MA ITER profile is much lower compared to that of Luxon and Brown and Taylor State demonstrating that the initial solution guess is a decent approximation for the 15MA ITER case.
Finally, the computed solutions corresponding to the Luxon and Brown model, Taylor state equilibrium, and 15MA ITER case
are shown in Figures~\ref{fig:lb}, \ref{fig:taylor}, and \ref{fig:empirical}.

In summary, we found that the proposed block triangular preconditioner with 3 W-cycles performs
well in most cases,  providing a good balance between 
the number of \das{inner AMG iterations and outer FGMRES iterations}.
In the rare cases when the W-cycle preconditioner fails,
a block triangular preconditioner with 3 to 5 V-cycles can be a more robust choice.
We also found that fewer nonlinear iterations are required after each application of AMR.

\subsection{Adaptive Mesh Refinement Case Studies}

\begin{figure}[htb]
  \centering
  \includegraphics[width=1\textwidth, trim={2cm 1.5cm 2.75cm 0cm}, clip]{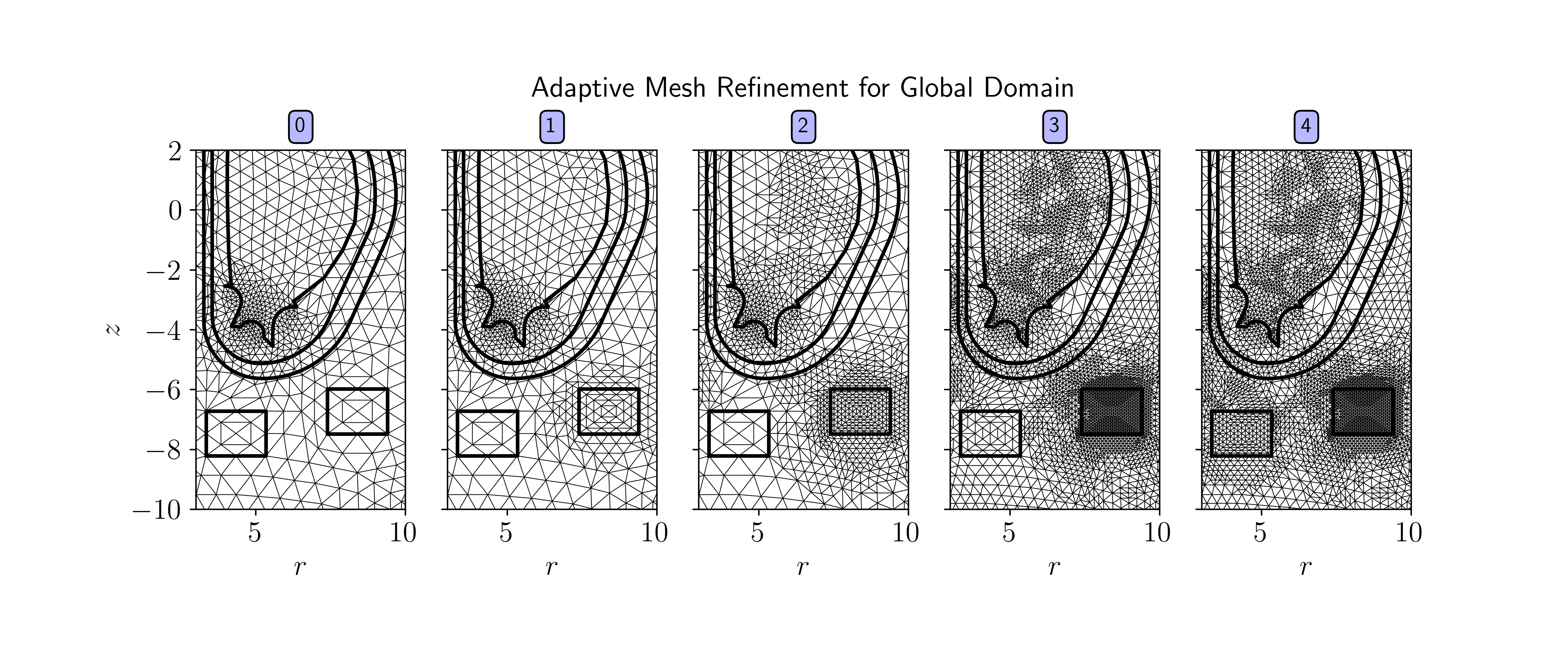}
  \caption{Four levels of AMR applied globally to $\psi$ for the Taylor equilibrium solution.\label{fig:amrglobal}}
\end{figure}

We demonstrate two different strategies for adaptive mesh refinement.
The first approach is to apply the error estimator to the solution $\psi$ at the end of each
Newton iteration.
Based on our experimentation, the error estimator typically prioritizes the regions near the coils
due to the presence of sharper gradients in $\psi$.
Due to this, we set a more aggressive threshold according to Section~\ref{sec:amrdesc}
in the limiter region to balance refinement near the plasma and near the coils.
For the first example, we consider four levels of AMR for the Taylor state equilibrium.
The meshes are shown in Figure~\ref{fig:amrglobal}.
The refinement targeted regions of large gradients near the coils and specific patches in the
limiter region.

\begin{figure}[htb]
  \centering
  \includegraphics[width=1\textwidth, trim={2cm 1.5cm 2.75cm 0cm}, clip]{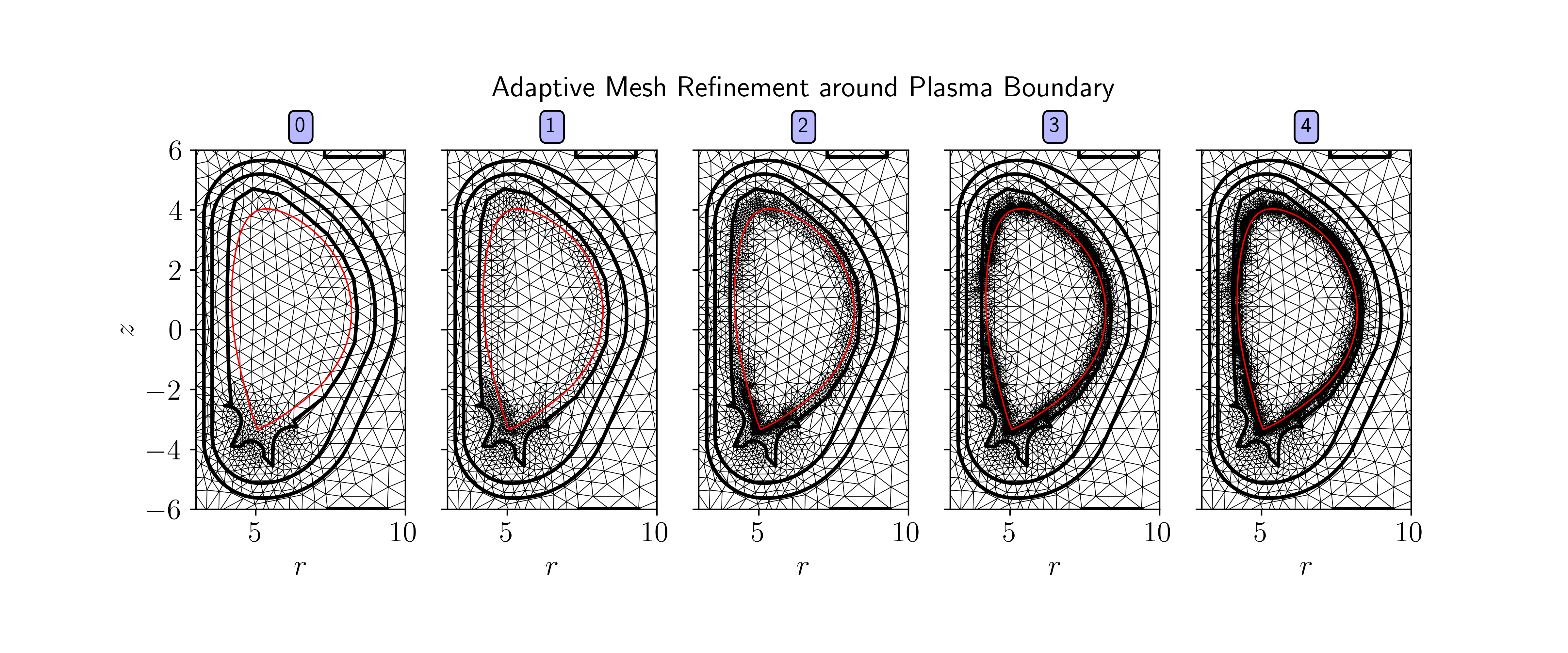}
    \caption{Four levels of AMR applied to $f(\psi)$, which is discontinuous across the plasma boundary. The solid red curve connects the control points used to optimize the plasma boundary.\label{fig:amrplasma}}
\end{figure}

We also demonstrate an alternative approach where we use $f(\psi)$, which is discontinuous across the
plasma boundary, in the error estimator.
This approach can be used to increase the resolution along the plasma boundary.
For this case, we also consider four levels of AMR for the Taylor state equilibrium.
The meshes are shown in Figure~\ref{fig:amrplasma}.
The refinement targets the region in the neighborhood around the control points.

\section{Conclusions}
\label{sec:conclusions}

In this work, we develop a Newton-based free-boundary Grad--Shafranov
solver based on adaptive finite elements, block preconditioning
strategies and algebraic multigrid.  A challenging optimization
problem to seek a free-boundary equilibrium has been successfully
addressed.  To seek an efficient and practical solver, we focus on
exploring different preconditioning strategies and leveraging advanced
algorithms such as AMR and scalable finite elements.  We identify
several simple but effective preconditioning options for the
linearized system.  We also demonstrate that the AMR algorithm helps
to resolve the local features and accelerates the nonlinear solver.
The resulting solver is particularly effective to seek the Taylor
state equilibrium that is needed in practical tokamak simulations of
vertical displacement events in a post-thermal-quench plasma.  Such an
equilibrium was impossible to solve by our previous Picard-based
solver.
In our Newton based approach, we demonstrated that the nonlinear residual
can be robustly reduced to 1e-6 and lower in a small handful of
iterations.
Although a good preconditioner strategy is identified, we note that
the performance of the preconditioner have some room for improvement.
A potential improvement might be to design a good smoother strategy
for the coupled system or its Schur complement.  Another future
direction is to incorporate the mesh generator into the solver loop,
which may address the potentially large change of the plasma domain
and can resolve the separatrix dynamically.  Such a feature has been
under active development in MFEM.

\section*{Acknowledgments}
DAS and QT would like to thank Veselin Dobrev for the advice on
implementing the double surface integral term.  This research used
resources provided by the Los Alamos National Laboratory Institutional
Computing Program, which is supported by the U.S. Department of
Energy's National Nuclear Security Administration under Contract
No.~89233218CNA000001, and the National Energy Research Scientific
Computing Center (NERSC), a U.S. Department of Energy Office of
Science User Facility located at Lawrence Berkeley National
Laboratory, operated under Contract No.~DE-AC02-05CH11231 using NERSC
award FES-ERCAP0028152.

\appendix
\section{Appendix}

We provide the details of computing the Gateaux semiderivatives of two exemplar terms using shape calculus.
The idea of using shape calculus in the Grad--Shafranov solver was first proposed in \cite{heumann2015quasi}.
However, the details of the derivation for the shape calculus were not provided in  \cite{heumann2015quasi}.
Here we provide all the necessary details and derivations for the purpose of reproduction.
\subsection{Gateaux Semiderivative of $a_{\rm plasma}$}
\label{sec:aplasma}
We first apply the definition of the Gateaux semiderivative,
\begin{align*}
  d_\psi(\varphi) [a_{\rm plasma}(\psi, v)] = -\frac{d}{d\epsilon}\bigg{(} \int_{\Omega_{\rm p}(\psi + \epsilon \varphi)} \bigg{(}& r p'(\psi + \epsilon \varphi) \\
  &\; + \frac{1}{\mu r} f(\psi + \epsilon \varphi) f'(\psi + \epsilon \varphi) \bigg{)} v \; dr dz, \bigg{)} \bigg{|}_{\epsilon=0}
\end{align*}
Using the Leibniz integral rule, the term becomes
\begin{align*}
  d_\psi(\varphi) [a_{\rm plasma}(\psi, v)]
  =\;&  - \int_{\Omega_{\rm p}(\psi)} d_{\psi}(\varphi)\left[ \left( r p'(\psi) + \frac{1}{\mu r} f(\psi) f'(\psi) \right) \right] v \; dr dz\\
     &- \int_{\partial \Omega_{\rm p}(\psi + \epsilon \varphi)} \bigg{(} r p'(\psi + \epsilon \varphi) \\
     &\qquad\qquad\qquad\quad+ \frac{1}{\mu r} f(\psi + \epsilon \varphi) f'(\psi + \epsilon \varphi) \bigg{)} \nv(\xv(\epsilon)) \cdot \frac{d \xv(\epsilon)}{d\epsilon} \, v \; d s \bigg{|}_{\epsilon=0},
\end{align*}
where $\nv(\xv(\epsilon))$ is the unit normal and at point $\xv(\epsilon)\in \partial \Omega(\psi+\epsilon\varphi)$.
Let $\xv_{\rm x}(\epsilon)$ be the saddle point of the field $\psi + \epsilon \varphi$.
The plasma boundary is implicitly defined by the equation
\begin{align*}
  \psi(\xv(\epsilon)) + \epsilon \varphi(\xv(\epsilon)) = \psi(\xv_{\rm x}(\epsilon)) + \epsilon \varphi(\xv_{\rm x}(\epsilon))
\end{align*}
By taking $\frac{d}{d\epsilon}$ and rearranging terms we have
\begin{align*}
  \left(\grad \psi(\xv(\epsilon)) + \epsilon \grad \varphi(\xv(\epsilon)) \right)\cdot \frac{d \xv(\epsilon)}{d \epsilon}
  &= \varphi(\xv_{\rm x}(\epsilon)) - \varphi
\end{align*}
Here, we used the fact that
$\grad \psi(\xv_{\rm x}(\epsilon)) + \epsilon \grad \varphi(\xv_{\rm x}(\epsilon)) = 0$
since, by definition, $\xv_{\rm x}(\epsilon)$ is a saddle point.
Note that the normal is given by
\begin{align*}
  \nv(\xv(\epsilon)) = \frac{\grad \psi(\xv(\epsilon)) + \epsilon \grad \varphi(\xv(\epsilon))}
  {|\grad \psi(\xv(\epsilon)) + \epsilon \grad \varphi(\xv(\epsilon))|}.
\end{align*}
Therefore,
\begin{align*}
  \nv(\xv(\epsilon)) \cdot \frac{d \xv(\epsilon)}{d\epsilon}
  = \frac{\varphi(\xv_{\rm x}(\epsilon)) - \varphi(\xv(\epsilon))}{|\grad \psi(\xv(\epsilon)) + \epsilon \grad \varphi(\xv(\epsilon))|},
\end{align*}
and
\begin{align*}
  d_\psi(\varphi) [a_{\rm plasma}(\psi, v)]
  =\;&  - \int_{\Omega_{\rm p}(\psi)} d_{\psi}(\varphi)\left[ \left( r p'(\psi) + \frac{1}{\mu r} f(\psi) f'(\psi) \right) \right] v \; dr dz\\
     &- \int_{\partial \Omega_{\rm p}(\psi)} \bigg{(} r p'(\psi) + \frac{1}{\mu r} f(\psi) f'(\psi) \bigg{)}
       |\grad \psi |^{-1} (\varphi(r_{\rm x}(\psi), z_{\rm x}(\psi)) - \varphi) \, v \; d s.
\end{align*}

\subsection{Gateaux Semiderivatives of $\psi_N(\psi)$}
\label{sec:semiempirical}
First, we use the quotient rule on~\eqref{eq:psiN}.
\begin{align}
  d_\psi(\varphi) [\psi_N(\psi)] =
  \frac{\left(\varphi-d_\psi(\varphi)[\psi_{\rm ma}]\right)(\psi_{\rm x} - \psi_{\rm ma})
  - (\psi - \psi_{\rm ma})(d_\psi(\varphi)[\psi_{\rm x}]
  - d_\psi(\varphi) [\psi_{\rm ma}])}{(\psi_{\rm x} - \psi_{\rm ma})^2}
  \label{eq:quotient}
\end{align}
Let $\xv_{\rm ma}(\epsilon)$ represent the location of the maximum point of the function $\psi+\epsilon \varphi$.
\begin{align*}
  d_\psi(\varphi) [\psi_{\rm ma}(\psi)] &= \frac{d}{d\epsilon}\left( \psi(\xv_{\rm ma}(\epsilon)) + \epsilon \varphi(\xv_{\rm ma}(\epsilon)\right)\bigg{|}_{\epsilon=0} \\
  &= \left(\grad \psi(\xv_{\rm ma}(\epsilon)) + \epsilon \grad\varphi(\xv_{\rm ma}(\epsilon)\right)\cdot \frac{d \xv(\epsilon)}{d \epsilon}\bigg{|}_{\epsilon=0} + \varphi(\xv_{\rm ma}(\epsilon))\bigg{|}_{\epsilon=0}
\end{align*}
Note that since $\xv_{\rm ma}$ corresponds to an extremum, the first term vanishes. Therefore,
\begin{align*}
  d_\psi(\varphi) [\psi_{\rm ma}(\psi, \varphi)] &= \varphi(\xv_{\rm ma}(0))
\end{align*}
A similar analysis follows for $d_\psi(\varphi) \psi_{\rm x}(\psi)$, since the first term vanishes at a
saddle point.
Therefore,
\begin{equation*}
  d_\psi(\varphi) [\psi_{\rm ma}(\psi)] = \varphi(r_{\rm ma}(\psi), z_{\rm ma}(\psi)),
  \qquad
  d_\psi(\varphi) [\psi_{\rm x}(\psi)] = \varphi(r_{\rm x}(\psi), z_{\rm x}(\psi)).
\end{equation*}
Substituting this into~\eqref{eq:quotient} and simplifying, we get
\begin{equation*}
  d_\psi \psi_N(\psi, \varphi) = \frac{\varphi-   (1 - \psi_N(\psi)) \varphi(r_{\rm ma}(\psi), z_{\rm ma}(\psi)) - \psi_N(\psi)\varphi(r_{\rm x}(\psi), z_{\rm x}(\psi))}{\psi_{\rm x} - \psi_{\rm ma}}.
\end{equation*}

\bibliographystyle{siamplain}
\bibliography{references}

\end{document}